\documentclass[12pt,reqno]{amsart}
\usepackage{latexsym, amsfonts, amsmath, amsthm, amssymb,amscd,color}
\usepackage[dvips]{graphicx}
\usepackage{eepic}
\usepackage{color}
\newtheorem{thm}{Theorem}[section]

\newtheorem{prop}[thm]{Proposition}

\newtheorem{cor}[thm]{Corollary}

\newtheorem{lem}[thm]{Lemma}

\numberwithin{equation}{section}
\def\pf{\noindent {\it Proof.} }

\def\1{$\bf{1}$}
\def\0{$\bf{0}$}
\def\c{{\rm{c}}}

\def\Rev{\rm{Rev}}
\def\rev{\rm{rev}}
\def\cro{\rm{cros}}
\def\ne{\rm{nest}}
\def\nec{\rm{ne}}
\def\sec{\rm{se}}
\def\luc{\rm{luc}}
\def\ruc{\rm{ruc}}
\def\le{\rm{left}}
\def\ri{\rm{right}}

\def\int{\rm{int}}

\def\Min{\rm{Min}}

\def\EC{\rm{EC}}
\def\ER{\rm{ER}}

\def\N{\mathcal{N}}
\def\G{\mathcal{G}}
\def\O{\mathcal{O}}
\def\C{\mathcal{C}}

\def\M{\mathcal{M}}
\def\P{\mathcal{P}}
\def\LP{\mathcal{LP}}

\title{Increasing and Decreasing sequences of length two in 01-fillings of moon polyominoes}
\author{Anisse Kasraoui}

\begin{document}
\maketitle \centerline{\small Universit\'e de Lyon ; \small
Universit\'e Lyon 1} \centerline{\small Institut Camille Jordan CNRS
UMR 5208} \centerline{\small 43, boulevard du 11 novembre 1918}
\centerline{\small F-69622, Villeurbanne Cedex }
\centerline{\small\texttt{anisse@math.univ-lyon1.fr}}

\begin{abstract}
We put recent results on the symmetry of the joint distribution of
the numbers of crossings and nestings of two edges over matchings
and set partitions in the larger context of the enumeration of
increasing and decreasing sequences of length $2$ in fillings of
moon polyominoes.
\end{abstract}

%%%%****%%%%****%%%%****%%%%****%%%%****%%%%****%%%%****%%%%*
%%%%
%%%% New Section
%%%%
%%%%****%%%%****%%%%****%%%%****%%%%****%%%%****%%%%****%%%%*

 \section{Introduction}

  The main purpose of this paper is to put recent results of Klazar and Noy~\cite{KlNo},
  Kasraoui and  Zeng~\cite{KaZe},
 and Chen, Wu and Yan~\cite{Chen1}, on the enumeration
of $2$-crossings and $2$-nestings in matchings, set partitions and
linked partitions in the larger context of enumeration of increasing
and decreasing chains in fillings of arrangements of cells.
 Our work is motivated by the recent paper of Krattenthaler~\cite{Kr} in which
results of Chen et al.~\cite{Chen2} on the symmetry of the crossing
number and nesting number in matchings and set partitions have been
extended in
a such context.\\[-0.3cm]

Let $G$ be a \emph{simple graph} (no multiple edges and loops) on
$[n]:=\{1,2,\ldots,n\}$. A graph will be represented by its set of
edges where the edge $\{i,j\}$ is written $(i,j)$ if $i<j$.

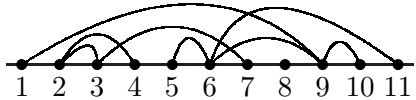
\begin{figure}[h!]
\begin{center}
{\setlength{\unitlength}{1mm}
\begin{picture}(50,10)(0,0)
\put(-2,0){\line(1,0){54}}
\put(0,0){\circle*{1,3}}\put(0,0){\makebox(0,-6)[c]{\small $1$}}
\put(5,0){\circle*{1,3}}\put(5,0){\makebox(0,-6)[c]{\small $2$}}
\put(10,0){\circle*{1,3}}\put(10,0){\makebox(0,-6)[c]{\small $3$}}
\put(15,0){\circle*{1,3}}\put(15,0){\makebox(0,-6)[c]{\small $4$}}
\put(20,0){\circle*{1,3}}\put(20,0){\makebox(0,-6)[c]{\small $5$}}
\put(25,0){\circle*{1,3}}\put(25,0){\makebox(0,-6)[c]{\small $6$}}
\put(30,0){\circle*{1,3}}\put(30,0){\makebox(0,-6)[c]{\small $7$}}
\put(35,0){\circle*{1,3}}\put(35,0){\makebox(0,-6)[c]{\small $8$}}
\put(40,0){\circle*{1,3}}\put(40,0){\makebox(0,-6)[c]{\small $9$}}
\put(45,0){\circle*{1,3}}\put(45,0){\makebox(0,-6)[c]{\small $10$}}
\put(50,0){\circle*{1,3}}\put(50,0){\makebox(0,-6)[c]{\small $11$}}
\qbezier(0,0)(25,16)(40,0)
\qbezier(5,0)(10,5)(10,0)\qbezier(20,0)(22,7)(25,0)
\qbezier(10,0)(20,10)(30,0)\qbezier(25,0)(30,15)(50,0)
\qbezier(25,0)(32.5,7)(40,0)\qbezier(40,0)(42,6)(45,0)
\qbezier(5,0)(10,8)(15,0)
\end{picture}}
\end{center}
\caption{The  graph
$\{(1,9),(2,3),(2,4),(3,7),(5,6),(6,9),(6,11),(9,10)\}$}\label{fig:graph}
\end{figure}
A sequence $(i_1,j_1),(i_2,j_2),\ldots,(i_k,j_k)$ of edges of $G$ is
said to be a \emph{$k$-crossing} if
$i_1<i_2<\cdots<i_k<j_1<j_2<\cdots<j_k$ and a \emph{$k$-nesting} if
$i_1<i_2<\cdots<i_k<j_k<\cdots<j_2<j_1$. If we draw the vertices of
$G$ in increasing order on a line and draw the arcs above the line
(see Figure~\ref{fig:graph} for an illustration), $k$-crossings and
$k$-nestings have a nice geometrical meaning. The largest $k$ for
which a graph $G$ has a $k$-crossing (resp., a $k$-nesting) is
denoted ${\cro}(G)$ (resp., ${\ne}(G)$) and called~\cite{Chen2} the
\textit{crossing number} (resp., \textit{nesting number}) of $G$.
The number of $k$-crossings (resp., $k$-nestings) of $G$ will be
denoted by ${\cro}_k(G)$ (resp., ${\ne}_k(G)$). A graph with no
$k$-crossing is called \emph{$k$-noncrossing} and a graph with no
$k$-nesting is called \emph{$k$-nonnesting}. As usual, a
$2$-noncrossing (resp., $2$-nonnesting) graph is just said to be
\emph{noncrossing} (resp., \emph{nonnesting}). Recently, there has
been an increasing interest in studying crossings and nestings in
matchings, set partitions, linked partitions and permutations (see
e.g. \cite{Bou,Chen1,Chen2,Co,Mi1,Mi2,KaZe,KlNo,Poz}).

\begin{figure}[h!]
\begin{center}
{\setlength{\unitlength}{1mm}
\begin{picture}(50,10)(0,0)
\put(-2,0){\line(1,0){54}}
\put(0,0){\circle*{1,3}}\put(0,0){\makebox(0,-6)[c]{\small 1}}
\put(5,0){\circle*{1,3}}\put(5,0){\makebox(0,-6)[c]{\small 2}}
\put(10,0){\circle*{1,3}}\put(10,0){\makebox(0,-6)[c]{\small 3}}
\put(15,0){\circle*{1,3}}\put(15,0){\makebox(0,-6)[c]{\small 4}}
\put(20,0){\circle*{1,3}}\put(20,0){\makebox(0,-6)[c]{\small 5}}
\put(25,0){\circle*{1,3}}\put(25,0){\makebox(0,-6)[c]{\small 6}}
\put(30,0){\circle*{1,3}}\put(30,0){\makebox(0,-6)[c]{\small 7}}
\put(35,0){\circle*{1,3}}\put(35,0){\makebox(0,-6)[c]{\small 8}}
\put(40,0){\circle*{1,3}}\put(40,0){\makebox(0,-6)[c]{\small 9}}
\put(45,0){\circle*{1,3}}\put(45,0){\makebox(0,-6)[c]{\small 10}}
\put(50,0){\circle*{1,3}}\put(50,0){\makebox(0,-6)[c]{\small 11}}
\qbezier(0,0)(20,12)(40,0) \qbezier(40,0)(42,5)(45,0)
\qbezier(5,0)(7,5)(10,0)\qbezier(20,0)(22,5)(25,0)
\qbezier(10,0)(20,8)(30,0)\qbezier(25,0)(37,10)(50,0)
\end{picture}
}
\end{center}
\caption{Standard representation of
$\pi=\{1,9,10\}\{2,3,7\}\{4\}\{5,6,11\}\{8\}$}\label{fig:partition-matching}
\end{figure}
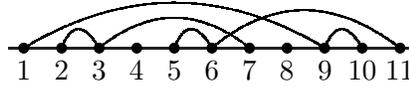

A (\emph{set}) \emph{partition} of $[n]$ is a collection of
non-empty pairwise disjoint sets, called blocks, whose union is
$[n]$. A (\emph{complete}) \emph{matching} of $[n]$ is just a set
partition whose each block contains exactly two elements. The set of
all set partitions and matchings of~$[n]$ will be denoted
respectively by $\P_n$ and $\M_n$. Set partitions (and thus
matchings) have a natural graphical representation, called
\emph{standard representation}. To each set partition $\pi$ of
$[n]$, one associates the graph $St_{\pi}$ on $[n]$ whose edge set
consists of arcs joining the elements of each block in numerical
order. For instance, the standard representation of the set
partition $\pi=\{\{1,9,10\},\{2,3,7\},\{4\},\{5,6,11\},\{8\}\}$ is
the graph on $\{1,2,\ldots,11\}$
$St_{\pi}=\{(1,9),(9,10),(2,3),(3,7), (5,6), (6,11)\}$ drawn in
Figure~\ref{fig:partition-matching}.

Throughout this paper, set partitions (and matchings) will be
identified with their standard representation. It is well-known that
the number of noncrossing matchings of $[2n]$ equals the number of
nonnesting matchings of $[2n]$, and that the number of noncrossing
partitions of $[n]$ equals the number of nonnesting partitions of
$[n]$ (and these are the $n$-th Catalan number), i.e.
\begin{align}
|\{M\in\M_{2n} :{\cro}_2(M)=0\}|&=|\{M\in\M_{2n} :{\ne}_2(M)=0\}|,\label{eq:noncrossing-match}\\
|\{\pi\in\P_{n} :{\cro}_2(\pi)=0\}|&=|\{\pi\in\P_{n}
:{\ne}_2(\pi)=0\}|\label{eq:noncrossing-part}.
\end{align}

In recent works, two generalizations of the latter identities have
been investigated. The first one is an extension of results obtained
by De Sainte-Catherine, and Klazar and Noy on the distributions of
$2$-crossings and $2$-nestings on matchings. In her
thesis~\cite{Sai}, De Sainte-Catherine have shown that the
statistics ${\cro}_2$ and ${\ne}_2$ are equidistributed over all
matchings of $[2n]$, that is for any integer $\ell\geq0$,
\begin{equation}
|\{M\in\M_{2n} :{\cro}_2(M)=\ell\}|=|\{M\in\M_{2n}
:{\ne}_2(M)=\ell\}|,
\end{equation}
i.e. in other words,
\begin{equation}\label{eq:equi-match}
\sum_{M\in\M_{2n}}p^{{\cro}_2(M)}=\sum_{M\in\M_{2n}}p^{{\ne}_2(M)}.
\end{equation}
 Klazar and Noy~\cite{KlNo} have shown that \eqref{eq:noncrossing-match} is even
 more true because the distribution of the joint
 statistic $({\cro}_2,{\ne}_2)$ is symmetric over $\M_{2n}$ that is
\begin{equation}\label{eq:sym-match}
\sum_{M\in\M_{2n}}p^{{\cro}_2(M)}q^{{\ne}_2(M)}=\sum_{M\in\M_{2n}}p^{{\ne}_2(M)}q^{{\cro}_2(M)}.
\end{equation}

Equations \eqref{eq:noncrossing-part} and
\eqref{eq:equi-match}\eqref{eq:sym-match} motivates Kasraoui and
Zeng to pose and solve the following questions: Are the statistics
${\cro}_2$ and ${\ne}_2$ equidistributed over all partitions of
$[n]$? Is the distribution of the joint statistic
$({\cro}_2,{\ne}_2)$ symmetric over all partitions of $[n]$? For
$S,T$ two subsets of $[n]$, let $\P_n(S,T)$ be the set of all
partitions of $[n]$ whose the set of lefthand (resp., righthand)
endpoints of the arcs of $\pi$ is equal to $S$ (resp., $T$). For
instance, the set partition drawn in
Figure~\ref{fig:partition-matching} belong to $\P_n(S,T)$, with
$S=\{1,2,3,5,6,9\}$ and $T=\{3,6,7,9,10,11\}$. Generalizing Klazar
and Noy's result \eqref{eq:sym-match}, Kasraoui and Zeng~\cite{KaZe}
have proved that the distribution of the joint statistic
$({\cro}_2,{\ne}_2)$ is symmetric over each $\P_n(S,T)$ (and thus,
over $\P_n$ and $\M_{n}$), that~is
\begin{equation}\label{eq:sym-partition}
\sum_{\pi\in\P_{n}(S,T)}p^{{\cro}_2(\pi)}q^{{\ne}_2(\pi)}
=\sum_{\pi\in\P_{n}(S,T)}p^{{\ne}_2(\pi)}q^{{\cro}_2(\pi)}.
\end{equation}
Note that recently, Chen, Wu and Yan~\cite{Chen1} have generalized
 the above result (although it is not explicitly stated) by
considering linked set partitions (see \eqref{eq:sym-linked}).

The second generalization of \eqref{eq:equi-match} and
\eqref{eq:sym-match} is due to  Chen, Deng, Du, Stanley and
Yan~\cite{Chen2}. It states, remarkably, that for any $k\geq2$ the
number of $k$-noncrossing partitions (resp. matchings) of $[n]$
equals the number of $k$-nonnesting partitions (resp. matchings)
of~$[n]$. More generally, Chen et al.~\cite{Chen2} have proved that
the distribution of the joint statistic $({\cro},{\ne})$ is
symmetric over each $\P_n(S,T)$, that is
\begin{equation}\label{eq:sym-partition-chen et al}
\sum_{\pi\in\P_{n}(S,T)}p^{{\cro}(\pi)}q^{{\ne}(\pi)}
=\sum_{\pi\in\P_{n}(S,T)}p^{{\ne}(\pi)}q^{{\cro}(\pi)}.
\end{equation}

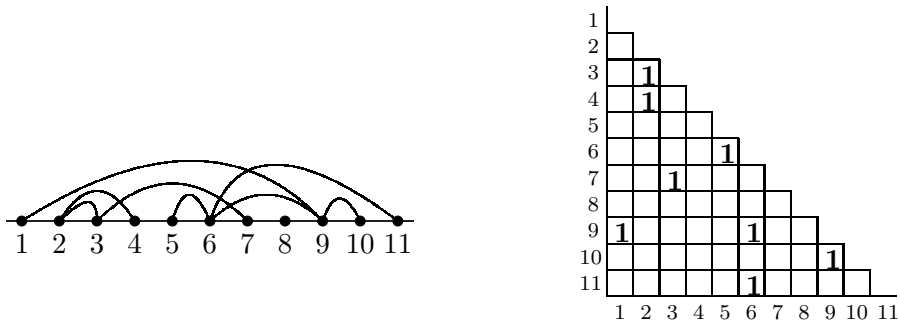
\begin{figure}[h!]
\begin{center}
{\setlength{\unitlength}{1mm}
\begin{picture}(50,30)(0,-10)
\put(-2,0){\line(1,0){54}}
\put(0,0){\circle*{1,3}}\put(0,0){\makebox(0,-6)[c]{\small $1$}}
\put(5,0){\circle*{1,3}}\put(5,0){\makebox(0,-6)[c]{\small $2$}}
\put(10,0){\circle*{1,3}}\put(10,0){\makebox(0,-6)[c]{\small $3$}}
\put(15,0){\circle*{1,3}}\put(15,0){\makebox(0,-6)[c]{\small $4$}}
\put(20,0){\circle*{1,3}}\put(20,0){\makebox(0,-6)[c]{\small $5$}}
\put(25,0){\circle*{1,3}}\put(25,0){\makebox(0,-6)[c]{\small $6$}}
\put(30,0){\circle*{1,3}}\put(30,0){\makebox(0,-6)[c]{\small $7$}}
\put(35,0){\circle*{1,3}}\put(35,0){\makebox(0,-6)[c]{\small $8$}}
\put(40,0){\circle*{1,3}}\put(40,0){\makebox(0,-6)[c]{\small $9$}}
\put(45,0){\circle*{1,3}}\put(45,0){\makebox(0,-6)[c]{\small $10$}}
\put(50,0){\circle*{1,3}}\put(50,0){\makebox(0,-6)[c]{\small $11$}}
\qbezier(0,0)(25,16)(40,0)
\qbezier(5,0)(10,5)(10,0)\qbezier(20,0)(22,7)(25,0)
\qbezier(10,0)(20,10)(30,0)\qbezier(25,0)(30,15)(50,0)
\qbezier(25,0)(32.5,7)(40,0)\qbezier(40,0)(42,6)(45,0)
\qbezier(5,0)(10,8)(15,0)
\end{picture}}\hspace{2.5cm}
{\setlength{\unitlength}{0.7mm}
\begin{picture}(50,50)(0,0)
%%%%%%%%%%%%%%%%%%%%%%%%%%%%%%%%%%%%%%%%DIAGRAMME
\put(0,0){\line(1,0){55}}\put(0,5){\line(1,0){50}}\put(0,10){\line(1,0){45}}
\put(0,15){\line(1,0){40}}\put(0,20){\line(1,0){35}}\put(0,25){\line(1,0){30}}
\put(0,30){\line(1,0){25}}\put(0,35){\line(1,0){20}}\put(0,40){\line(1,0){15}}
\put(0,45){\line(1,0){10}}\put(0,50){\line(1,0){5}}
\put(0,0){\line(0,1){55}}\put(5,0){\line(0,1){50}}\put(10,0){\line(0,1){45}}
\put(15,0){\line(0,1){40}}\put(20,0){\line(0,1){35}}\put(25,0){\line(0,1){30}}
\put(30,0){\line(0,1){25}}\put(35,0){\line(0,1){20}}\put(40,0){\line(0,1){15}}
\put(45,0){\line(0,1){10}}\put(50,0){\line(0,1){5}}
%%%%%%%%%%%%%%%%%%%%%%%%%%%%%%%%%%%%%%%%REMPLISSAGE
\put(0,10){\makebox(6,4)[c]{\small \1}}
\put(5,35){\makebox(6,4)[c]{\small \1}}
\put(5,40){\makebox(6,4)[c]{\small \1}}
\put(10,20){\makebox(6,4)[c]{\small \1}}
\put(20,25){\makebox(6,4)[c]{\small \1}}
\put(25,0){\makebox(6,4)[c]{\small \1}}
\put(25,10){\makebox(6,4)[c]{\small \1}}
\put(40,5){\makebox(6,4)[c]{\small \1}}
%%%%%%%%%%%%%%%%%%%%%%%%%%%%%%%%%%%%%%%%ETIQUETAGE DES COLONNES ET LIGNES
\put(0,0){\makebox(5,-6)[c]{\tiny $1$}}
\put(5,0){\makebox(5,-6)[c]{\tiny $2$}}
\put(10,0){\makebox(5,-6)[c]{\tiny $3$}}
\put(15,0){\makebox(5,-6)[c]{\tiny $4$}}
\put(20,0){\makebox(5,-6)[c]{\tiny $5$}}
\put(25,0){\makebox(5,-6)[c]{\tiny $6$}}
\put(30,0){\makebox(5,-6)[c]{\tiny $7$}}
\put(35,0){\makebox(5,-6)[c]{\tiny $8$}}
\put(40,0){\makebox(5,-6)[c]{\tiny $9$}}
\put(45,0){\makebox(5,-6)[c]{\tiny $10$}}
\put(51,0){\makebox(5,-6)[c]{\tiny $11$}}
%%%%%%%%%%%%
\put(0,0){\makebox(-6,5)[c]{\tiny $11$}}
\put(0,5){\makebox(-6,5)[c]{\tiny $10$}}
\put(0,10){\makebox(-5,5)[c]{\tiny $9$}}
\put(0,15){\makebox(-5,5)[c]{\tiny $8$}}
\put(0,20){\makebox(-5,5)[c]{\tiny $7$}}
\put(0,25){\makebox(-5,5)[c]{\tiny $6$}}
\put(0,30){\makebox(-5,5)[c]{\tiny $5$}}
\put(0,35){\makebox(-5,5)[c]{\tiny $4$}}
\put(0,40){\makebox(-5,5)[c]{\tiny $3$}}
\put(0,45){\makebox(-5,5)[c]{\tiny $2$}}
\put(0,50){\makebox(-5,5)[c]{\tiny $1$}}
\end{picture}}
\end{center}
\caption{A graph and the corresponding
$01$-filling}\label{fig:part-filling}
\end{figure}

In the recent paper~\cite{Kr}, Krattenthaler have put Chen et al's
result~\eqref{eq:sym-partition-chen et al} in the larger context of
the enumeration of increasing and decreasing chains in fillings of
ferrers shapes. First recall the correspondence between simple
graphs of $[n]$ and $01$-fillings of $\Delta_n$, the triangular
shape with $n-1$ cells in the bottom row, $n-2$ cells in the row
above, etc., and $1$ cell in the top-most row. See
Figure~\ref{fig:part-filling} for an example in which $n =11$ (the
filling and labeling of the corners should be ignored at this point.
For convenience, we also joined pending edges at the right and at
the top of $\Delta_n$). Let $G$ be a simple graph of $[n]$. The
correspondence consists in labeling in increasing order columns from
left to right by $\{1,2,\ldots,n\}$ and rows from top to bottom by
$\{1,2,\ldots,n\}$. Then assign the value 1 to the cell on column
labeled~$i$ and row labeled~$j$ if and only if $(i,j)$ is an edge of
$G$. An illustration is given in Figure~\ref{fig:part-filling}.

It is obvious to see that in this correspondence a $k$-crossing
(resp., $k$-nesting) corresponds to a SE-chain (resp., NE-chain) of
length $k$ such that the smallest rectangle containing the chain is
contained in~$\Delta_n$, a SE-chain (resp., NE-chain) of length $k$
being a sequence of $k$ 1's in the filling such that any 1 in the
sequence is below and to the right (resp., above and to the right)
of the preceding 1 in the sequence. Moreover, it is obvious that
this correspondence establishes a bijection between set partitions
of $[n]$ and ${\N}(\Delta_n)$, the set of all $01$-fillings of
$\Delta_n$ in which every row and every column contains at most
one~1. Therefore, Kasraoui and Zeng's
result~\eqref{eq:sym-partition} and Chen et al's
result~\eqref{eq:sym-partition-chen et al} can be viewed as a
property of symmetry of NE-chains and SE-chains over
${\N}(\Delta_n)$. Given a $01$-filling $F$ of $\Delta_n$, denote by
${\sec}(F)$ (resp., ${\nec}(F)$) the maximal $k$ such that $F$ has a
SE-chain (resp., NE-chain) of length $k$, the smallest rectangle
containing the chain being contained in~$F$ and by ${\sec}_2(F)$
(resp., ${\nec}_2(F)$) the number of SE-chains (resp., NE-chains) of
length $2$ such that the smallest rectangle containing the chain is
contained in~$F$. Then the symmetry of the distributions of the
joint statistics $({\cro}_2,{\ne}_2)$ and $({\cro},{\ne})$ over set
partitions can be reformulated respectively as follows:
\begin{align}
\sum_{F\in\,{\N}(\Delta_n)}p^{{\nec}_2(F)}q^{{\sec}_2(F)}
&=\sum_{F\in\,{\N}(\Delta_n)}p^{{\sec}_2(F)}q^{{\nec}_2(F)}\label{eq:intro1}\\
\sum_{F\in\,{\N}(\Delta_n)}p^{{\nec}(F)}q^{{\sec}(F)}
&=\sum_{F\in\,{\N}(\Delta_n)}p^{{\sec}(F)}q^{{\nec}(F)}.\label{eq:intro2}
\end{align}

In the recent paper~\cite{Kr}, Krattenthaler have shown
that~\eqref{eq:intro2} remains true if we replace $\Delta_n$ by any
ferrers shapes and he proposed to investigate more general
arrangements. This was done successfully by Rubey~\cite{Ru} for moon
polyominoes. It is thus natural to ask if~\eqref{eq:intro1} remains
true when we replace $\Delta_n$ by any ferrers shape or, more
generally by any moon polyomino (We will answer this question by the
affirmative): this is the original motivation of our paper.

\section{The main results}

A \emph{polyomino} is an arrangement of square cells. It is
\emph{convex} if along any row of cells and along any column of
cells there is no hole. It is \emph{intersection free} if any two
rows are comparable, i.e., one row can be embedded in the other by
applying a vertical shift. Equivalently, it is \emph{intersection
free} if any two columns are comparable, i.e., one row can be
embedded in the other by applying an horizontal shift. A \emph{moon
polyomino} is a convex and intersection free polyomino. An
illustration is given in Figure~\ref{fig:moonpolyominoe-filling}.
\begin{center}
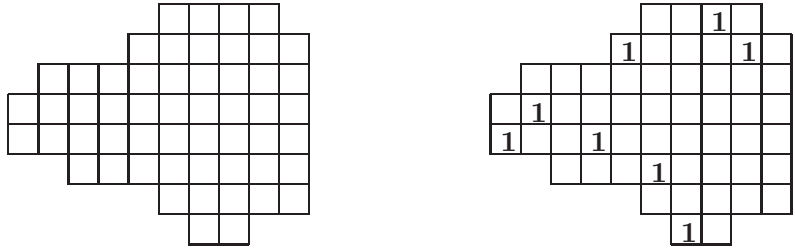
\begin{figure}[h!]
 {\setlength{\unitlength}{0.8mm}
\begin{picture}(50,50)(0,0)
%%%%%%%%%%%%%%%%%%%%%%%%%%%%%%%%%%%%%%%%DIAGRAMME
\put(30,0){\line(1,0){10}} \put(25,5){\line(1,0){25}}
\put(10,10){\line(1,0){40}} \put(0,15){\line(1,0){50}}
\put(0,20){\line(1,0){50}} \put(0,25){\line(1,0){50}}
\put(5,30){\line(1,0){45}} \put(20,35){\line(1,0){30}}
\put(15,30){\line(1,0){10}} \put(25,40){\line(1,0){20}}
\put(0,15){\line(0,1){10}}\put(5,15){\line(0,1){15}}
\put(10,10){\line(0,1){20}}\put(15,10){\line(0,1){20}}\put(20,10){\line(0,1){25}}
\put(25,5){\line(0,1){35}}\put(30,0){\line(0,1){40}}\put(35,0){\line(0,1){40}}\put(40,0){\line(0,1){40}}
\put(45,5){\line(0,1){35}} \put(50,5){\line(0,1){30}}
\end{picture}}
\hspace{2cm} {\setlength{\unitlength}{0.8mm}
\begin{picture}(50,50)(0,0)
%%%%%%%%%%%%%%%%%%%%%%%%%%%%%%%%%%%%%%%%coloriage
%\put(10,10){\color{red}{\rule{4mm}{16mm}}}
%\put(45,5){\color{red}{\rule{4mm}{24mm}}}
%%%%%%%%%%%%%%%%%%%%%%%%%%%%%%%%%%%%%%%%DIAGRAMME
\put(30,0){\line(1,0){10}} \put(25,5){\line(1,0){25}}
\put(10,10){\line(1,0){40}} \put(0,15){\line(1,0){50}}
\put(0,20){\line(1,0){50}} \put(0,25){\line(1,0){50}}
\put(5,30){\line(1,0){45}} \put(20,35){\line(1,0){30}}
\put(15,30){\line(1,0){10}} \put(25,40){\line(1,0){20}}
\put(0,15){\line(0,1){10}}\put(5,15){\line(0,1){15}}
\put(10,10){\line(0,1){20}}\put(15,10){\line(0,1){20}}\put(20,10){\line(0,1){25}}
\put(25,5){\line(0,1){35}}\put(30,0){\line(0,1){40}}\put(35,0){\line(0,1){40}}\put(40,0){\line(0,1){40}}
\put(45,5){\line(0,1){35}} \put(50,5){\line(0,1){30}}
%%%%%%%%%%%%%%%%%%%%%%%%%%%%%%%%%%%%%%%%REMPLISSAGE
\put(30,0){\makebox(6,4)[c]{\small \1}}
\put(25,10){\makebox(6,4)[c]{\small \1}}
\put(15,15){\makebox(6,4)[c]{\small\1}}
\put(0,15){\makebox(6,4)[c]{\small \1}}
\put(5,20){\makebox(6,4)[c]{\small\1}}
\put(20,30){\makebox(6,4)[c]{\small\1}}
\put(40,30){\makebox(6,4)[c]{\small \1}}
\put(35,35){\makebox(6,4)[c]{\small \1}}
\end{picture}}
\caption{A moon polyomino $T$ and a $01$-filling of
$T$.}\label{fig:moonpolyominoe-filling}
\end{figure}
\end{center}
 Let $T$ be a moon polyomino. A $01$-\emph{filling} $F$ of $T$ consists
of assigning 0 or 1 to each cell. For convenience, we will omit the
0's when we draw the fillings. See
Figure~\ref{fig:moonpolyominoe-filling}. The set of all
$01$-fillings of $T$ will be denoted $\N^{01}(T)$. Recall that a
SE-chain (resp., NE-chain) of length $k$ in a $01$-filling~$F$ of
$T$ is a sequence of $k$ 1's in the filling such that any 1 in the
sequence is below and to the right (resp., above and to the right)
of the preceding 1 in the sequence. A SE-chain (resp., NE-chain) of
length $2$ such that the smallest rectangle containing the chain is
contained in $F$ is said to be a \emph{descent} (resp., an
\emph{ascent}). We will denote by ${\sec}_2(F)$ and ${\nec}_2(F)$
the number of descents and ascents in~$F$. For instance, if $F$ is
the filling drawn in Figure~\ref{fig:part-filling}, we have
${\nec}_2(F)=6$ and ${\sec}_2(F)=4$, while for the filling in
Figure~\ref{fig:moonpolyominoe-filling} we have
${\nec}_2(F)={\sec}_2(F)=4$. It is natural in view of the results
presented in the introduction to ask if the statistics ${\cro}_{2}$
and ${\ne}_2$ are equidistributed over all simple graphs of $[n]$,
or equivalently, if the statistics ${\sec}_{2}$ and ${\nec}_2$ are
equidistributed over $\N^{01}(\Delta_n)$ for any nonnegative
integer~$n$. More generally, one can ask if  the statistics
${\sec}_{2}$ and ${\nec}_2$ are equidistributed over $\N^{01}(T)$
for any moon polyomino~$T$. The answer to these questions is no by
means of Proposition~\ref{prop:nonsym-graphe}. However, it appears
that for particular $01$-fillings we have such an equidistribution
and even more, namely the symmetry of the distribution of the joint
statistic $({\nec}_{2},{\sec}_2)$. Let $\N^c(T)$ (resp., $\N^r(T)$)
be the set of all $01$-fillings  of $T$ with at most one $1$ in each
column (resp., row), and $\N(T):=\N^c(T)\cap\N^r(T)$ the set of all
$01$-fillings of $T$ with at most one $1$ in each column and in each
row. Then our main result can be stated as follows.

\begin{thm}\label{thm:sym-filling}
For any moon polyomino $T$, the distribution of the joint statistic
$({\nec}_2,{\sec}_2)$ over any ${\Large
\mathcal{B}}\in\{\N(T),\N^c(T),\N^r(T)\}$ is symmetric, or
equivalently
\begin{align}
\sum_{F\in\,{\Large \mathcal{B}}}p^{{\nec}_2(F)}q^{{\sec}_2(F)}
&=\sum_{F\in\,{\Large \mathcal{B}}}p^{{\sec}_2(F)}q^{{\nec}_2(F)}.
\end{align}
\end{thm}

In fact we have obtained much stronger results. Before stating these
results, we need to introduce some definitions.
 Let $T$ be a moon
polyomino with $s$ rows and $t$ columns. By convention, we always
label the rows of $T$ from top to bottom in increasing order by
$\{1,2,\ldots,s\}$ and the columns of $T$ from left to right in
increasing order by $\{1,2,\ldots,t\}$. The row labeled~$i$ and the
column labeled~$j$ will be denoted respectively by $R_i$ and $C_j$.
The \emph{length-row sequence} of $T$, denoted $r(T)$, is the
sequence $(r_1,r_2,\ldots,r_s)$ where $r_i$ is the length (i.e., the
number of cells) of the row $R_i$. Similarly, the
\emph{length-column sequence} of $T$, denoted $c(T)$, is the
sequence $(c_1,c_2,\ldots,c_t)$ where $c_i$ is the length of the
column $C_i$. Clearly, the length-row and length-column sequence of
any moon polyomino are always unimodal sequences, that is there
exist (unique) integers $i_0$ and $j_0$ such that $r_1\leq
r_2\leq\cdots\leq r_{i_0}>r_{i_0+1}\geq\cdots\geq r_s$ and $c_1\leq
c_2\leq\cdots\leq c_{j_0}>c_{j_0+1}\geq\cdots\geq c_t$. The
\emph{upper part of~$T$}, denoted $Up(T)$, is the set of rows $R_i$
with $1\leq i\leq i_0$, and the \emph{lower part}, denoted $Low(T)$,
the set of rows $R_i$, $i_0+1\leq i\leq s$. Similarly, the
\emph{left part} of $T$, denoted $Left(T)$, is the set of columns
$C_i$ with $1\leq i\leq j_0$, and the \emph{right part}, denoted
$Right(T)$, the set of columns $C_i$, $j_0+1\leq i\leq t$. For
instance, if~$T$ is the moon polyomino in
Figure~\ref{fig:moonpolyominoe-filling}, we have
$r(T)~=~(4,6,9,10,10,8,5,2)$, $Up(T)=\{R_i: 1\leq i\leq 5\}$ and
$Low(T)=\{R_6,R_7,R_8\}$, and $c(T)=(2,3,4,4,5,7,8,8,7,6)$,
$Left(T)=\{C_i: 1\leq i\leq 8\}$ and $Right(T)=\{C_9,C_{10}\}$.
Define the relation $\prec$ on the rows of $T$ as follows: $R_i\prec
R_j$ if and only if
\begin{itemize}
\item $r_i<r_j$, or
\item $r_i=r_j$, $R_i\in Up(T)$ and $R_j\in Low(T)$, or
\item $r_i=r_j$, $R_i,R_j\in Up(T)$ and $R_i$ is above $R_j$, or
\item $r_i=r_j$, $R_i,R_j\in Low(T)$ and $R_i$ is below $R_j$.
\end{itemize}
Similarly, define the relation $\prec$ (for convenience, we use the
same symbol than for rows) on the columns of $T$ defined by
$C_i\prec C_j$ if and only if
\begin{itemize}
\item $c_i<c_j$, or
\item $c_i=c_j$, $C_i\in Left(T)$ and $C_j\in Right(T)$, or
\item $c_i=c_j$, $C_i,C_j\in Left(T)$ and $C_i$ is to the left of $C_j$, or
\item $c_i=c_j$, $C_i,C_j\in Right(T)$ and $C_i$ is to the right of $C_j$.
\end{itemize}
It is easy to check that the relation $\prec$ is a total order both
on rows and columns of $T$. For instance, if $T$ is the moon
polyomino in Figure~\ref{fig:moonpolyominoe-filling} we have
$R_8\prec R_1\prec R_7\prec R_2\prec R_6\prec R_3\prec R_4\prec R_5$
and  $C_1\prec C_2\prec
C_3\prec C_4\prec C_5\prec C_{10}\prec C_6\prec C_9\prec C_7\prec C_8$.\\[-0.3cm]

Let $F$ be a $01$-filling of $T$. A cell of $F$ is said to be
\emph{empty} if it has been assigned the value~$0$. We also say that
a row (resp., column) of $F$ is empty if all its cells are empty.
The indices of the empty rows and columns of $F$ are denoted
${\ER}(F)$ and ${\EC}(F)$, respectively. For instance if $F$ is the
$01$-filling given in Figure~\ref{fig:moonpolyominoe-filling}, then
${\ER}(F)=\{3,7\}$  and ${\EC}(F)=\{3,10\}$.

Given a $s$-uple ${\bf m}=(m_1,\ldots,m_s)$ of positive integers and
$A$ a subset of $[t]$, we denote by $\N^c(T,{\bf m})$ the set of
$01$-fillings in $\N^c(T)$ with exactly $m_i$ $1$'s in row $R_i$ and
by $\N^c(T,{\bf m};A)$ the set of fillings $F$ in $\N^c(T,{\bf m})$
such that ${\EC}(F)=A$. For instance, the filling $F$ given in
Figure~\ref{fig:moonpolyominoe-filling} belong to $\N^c(T,{\bf
m};A)$, with ${\bf m}=~(1,2,0,1,2,1,0,1)$ and $A=\{3,10\}$.
Similarly, given a $t$-uple ${\bf n}=(n_1,\ldots,n_t)$ of positive
integers and a subset $B$ of $[s]$, we denote by $\N^r(T,{\bf n})$
the set of $01$-fillings in $\N^r(T)$ with exactly $n_i$ $1$'s in
column $C_i$ and by $\N^r(T,{\bf n};B)$ the set of fillings $F$ in
$\N^r(T,{\bf n})$ such that ${\ER}(F)=B$.
 Also, for $A,B$ two subsets of $[t]$ and $[s]$ respectively, we denote by
$\N(T;A,B)$ the set of $01$-fillings in $\N(T)$ such that
${\EC}(F)=A$ and ${\ER}(F)=B$.

For positive integers $n$ and $k$, let ${n\brack k}_{p,q}$ be the
$p,q$-Gaussian coefficient defined by
$$
{n\brack k}_{p,q}=\left\{
                   \begin{array}{ll}
                     \frac{[n]_{p,q}!}{[k]_{p,q}!\;[n-k]_{p,q}!}, & \hbox{if $0\leq k\leq n$;} \\
                     0, & \hbox{otherwise.}
                   \end{array}
                 \right.
$$
where, as usual in $p,q$-theory, the $p,q$-integer $[r]_{p,q}$ is
given by
$$[r]_{p,q}:=\frac{p^i-q^i}{p-q}=(p^{i-1}+p^{i-2}q+\cdots+p^jq^{i-j-1}+\cdots+pq^{i-2}+q^{i-1}),$$
and the $p,q$-factorial $[r]_{p,q}!$  by
$[r]_{p,q}!:=\prod_{i=1}^r[i]_{p,q}$.\\

  Let $T$ be a moon polyomino with
$s$ rows and $t$ columns, ${\bf m}=(m_1,\ldots,m_s)$ a $s$-uple of
positive integers, ${\bf n}=(n_1,\ldots,n_t)$ a $t$-uple of positive
integers, and $A,B$  two subsets of $[t]$ and $[s]$ respectively.
Suppose $R_{i_1}\prec R_{i_2} \prec \cdots \prec R_{i_s}$ and
$C_{j_1}\prec C_{j_2} \prec \cdots \prec C_{j_t}$. Then for
$u\in[s]$ and $v\in[t]$, define $h_{i_u}$ and $h'_{j_v}$ by
\begin{align}
h_{i_u}&=r_{i_u}-(m_{i_1}+m_{i_2}+\cdots+m_{i_{u-1}})-a_{{i_u}},\label{eq:h_i}\\
h'_{j_v}&=c_{j_v}-(n_{j_1}+n_{j_2}+\cdots+n_{j_{v-1}})-b_{{j_v}},\label{eq:h'i}
\end{align}
where  $r_{i_u}$ is the length of the row $R_{i_u}$ and $a_{{i_u}}$
is the number of indices $k\in A$ such that the column $C_k$
intersect the row $R_{i_u}$, and $c_{j_v}$ is the length of the
column $C_{j_v}$ and $b_{{j_v}}$ is the number of indices $k\in B$
such that the row $R_k$ intersect the column $C_{j_v}$.

 The following result gives the distributions of  the joint statistic $({\nec}_2,{\sec}_2)$
over $\N^c(T,{\bf m};A)$ and  $\N^r(T,{\bf n};B)$.

\begin{thm}\label{thm:distribution}
For any moon polyomino $T$ with $s$ rows and $t$ columns, the
distributions of the joint statistic $({\nec}_2,{\sec}_2)$ over each
$\N^c(T,{\bf m};A)$ and $\N^r(T,{\bf n};B)$ are given by
\begin{align}
\sum_{F\in\,\N^c(T,{\bf
m};A)}p^{{\nec}_2(F)}q^{{\sec}_2(F)}&=\prod_{d=1}^s{h_d\brack
m_d}_{p,q},\label{eq:symfilling}\\
\sum_{F\in\,\N^r(T,{\bf
n};B)}p^{{\nec}_2(F)}q^{{\sec}_2(F)}&=\prod_{d=1}^t{h'_d\brack
n_d}_{p,q},\label{eq:symfilling'}
\end{align}
where $h_d$ and $h'_d$ are defined by \eqref{eq:h_i} and
\eqref{eq:h'i}.
\end{thm}

 As an immediate consequence (take ${\bf m}\in\{0,1\}^s$ or ${\bf  n}\in\{0,1\}^t$ in
Theorem~\ref{thm:distribution}), we obtain the following result.
\begin{cor}
For any moon polyomino $T$ with $s$ rows and $t$ columns, the
distribution of the joint statistic $({\nec}_2,{\sec}_2)$ over
$\N(T;A,B)$ is given by
\begin{align}
\sum_{F\in\,\N(T;A,B)}p^{{\nec}_2(F)}q^{{\sec}_2(F)}
&=\prod_{d\in[s]\setminus B}[h_d]_{p,q}=\prod_{d\in[t]\setminus
A}[h'_d]_{p,q},
\end{align}
where $h_d$ and $h'_d$ are defined by \eqref{eq:h_i} and
\eqref{eq:h'i}.
\end{cor}

It is not difficult to show that \eqref{eq:symfilling} and
\eqref{eq:symfilling'} are equivalent. Indeed, if $T$ is a moon
polyomino with $s$ rows and $t$ columns, then the arrangement of
cells $T^t$ obtained from $T$ by rotation about $90^{\circ}$ is also
a moon polyomino but with $t$ rows and $s$ columns. Moreover it is
obvious that the application $F\mapsto F^t$ which associates to each
filling $F$ of $T$ the filling $F^t$ of $T^t$ obtained from $F$ by
rotation about $90^{\circ}$ establishes a bijection from
$\N^c(T,{\bf m};A)$ onto $\N^r(T^t,{\bf m};A)$ which sends the joint
statistic $({\nec}_2,{\sec}_2)$ onto $({\sec}_2,{\nec}_2)$ for
any~${\bf m}$ and~$A$. See Figure~\ref{fig:rotation-filling} for an
illustration.

\begin{center}
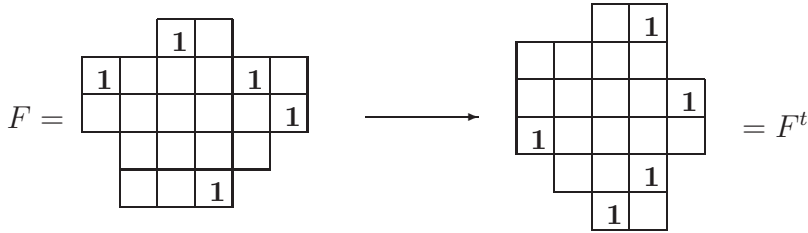
\begin{figure}[h!] {\setlength{\unitlength}{1mm}
\begin{picture}(35,25)(0,-3)
\put(-10,10){\makebox(6,4)[l]{$F=$}}
%%%%%%%%%%%%%%%%%%%%%%%%%%%%%%%%%%%%%%%%coloriage
%\put(5,0){\color{red}{\rule{5mm}{20mm}}}
%%%%%%%%%%%%%%%%%%%%%%%%%%%%%%%%%%%%%%%%DIAGRAMME
\put(5,0){\line(1,0){15}}\put(5,5){\line(1,0){20}}\put(0,10){\line(1,0){30}}
\put(0,15){\line(1,0){30}}\put(0,20){\line(1,0){30}}\put(10,25){\line(1,0){10}}
\put(0,10){\line(0,1){10}}\put(5,0){\line(0,1){20}}\put(10,0){\line(0,1){25}}
\put(15,0){\line(0,1){25}}\put(20,0){\line(0,1){25}}\put(25,5){\line(0,1){15}}
\put(30,10){\line(0,1){10}}
%%%%%%%%%%%%%%%%%%%%%%%%%%%%%%%%%%%%%%%%REMPLISSAGE
\put(10,20){\makebox(6,4)[c]{\small\1}}
\put(0,15){\makebox(6,4)[c]{\small\1}}
\put(20,15){\makebox(6,4)[c]{\small \1}}
\put(25,10){\makebox(6,4)[c]{\small \1}}
\put(15,0){\makebox(6,4)[c]{\small \1}}
\end{picture}}\hspace{2cm}
 {\setlength{\unitlength}{1mm}
\begin{picture}(35,35)(0,0)
\put(-20,15){\vector(1,0){15}}
 \put(30,12){\makebox(6,4)[l]{$=F^t$}}
%%%%%%%%%%%%%%%%%%%%%%%%%%%%%%%%%%%%%%%%coloriage
%\put(5,0){\color{red}{\rule{5mm}{20mm}}}
%%%%%%%%%%%%%%%%%%%%%%%%%%%%%%%%%%%%%%%%DIAGRAMME
\put(10,0){\line(1,0){10}}\put(5,5){\line(1,0){15}}\put(0,10){\line(1,0){25}}
\put(0,15){\line(1,0){25}}\put(0,20){\line(1,0){25}}\put(0,25){\line(1,0){20}}
\put(10,30){\line(1,0){10}}\put(0,10){\line(0,10){15}}\put(5,5){\line(0,1){20}}
\put(10,0){\line(0,1){30}}\put(15,0){\line(0,1){30}}\put(20,0){\line(0,1){30}}
\put(25,10){\line(0,1){10}}
%%%%%%%%%%%%%%%%%%%%%%%%%%%%%%%%%%%%%%%%REMPLISSAGE
\put(10,0){\makebox(6,4)[c]{\small\1}}
\put(15,5){\makebox(6,4)[c]{\small\1}}
\put(0,10){\makebox(6,4)[c]{\small \1}}
\put(20,15){\makebox(6,4)[c]{\small \1}}
\put(15,25){\makebox(6,4)[c]{\small \1}}
\end{picture}}
\caption{Rotation about $90^{\circ}$ for
filling.}\label{fig:rotation-filling}
\end{figure}
\end{center}
Since the $p,q$-integer $[n]_{p,q}$ is symmetric in the variables
$p$ and $q$ for any positive integer~$n$, we have the following
results.

\begin{cor}\label{thm:refinement main1}
For any moon polyomino $T$, the joint statistic
$({\nec}_2,{\sec}_2)$ is symmetrically distributed over each
$\N^c(T,{\bf m};A)$, $\N^r(T,{\bf n};B)$ and $\N(T;A,B)$. In
particular, it is symmetrically distributed over each $\N^c(T,{\bf
m})$ and $\N^r(T,{\bf n})$.
\end{cor}

For any positive integer $k$ denote by $\N^c(T;k)$, $\N^r(T;k)$ and
$\N(T;k)$ the set of $01$-fillings in $\N^c(T)$, $\N^r(T)$ and
$\N(T)$ with exactly $k$ ones, respectively.
\begin{cor}\label{thm:refinement main2}
For any moon polyomino $T$ and positive integer~$k$, the joint
statistic $({\nec}_2,{\sec}_2)$ is symmetrically distributed over
each $\N^c(T;k)$, $\N^r(T;k)$ and $\N(T;k)$. Summing over all
positive integers $k$, we recover Theorem~\ref{thm:sym-filling}.
\end{cor}

The paper is organized as follows. In Section~3, we show briefly how
the results on the symmetry of the joint statistic
$({\cro}_2,\ne_2)$ presented in the introduction can be obtained
from the above results. In section~4, we prove
Theorem~\ref{thm:sym-filling} and in Section~5, we present a
bijective proof of Corollary~\ref{thm:refinement main1}.
 Finally, we conclude this
paper with some remarks and problems.\\

We end this section by illustrating Theorem~\ref{thm:distribution}.
Suppose $T$ is the moon polyomino given below and $A=\{2\}$.
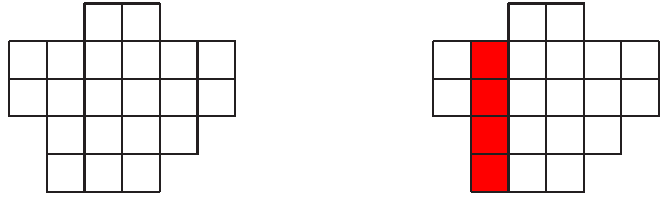
\begin{center}
 {\setlength{\unitlength}{1mm}
\begin{picture}(35,28)(0,0)
%%%%%%%%%%%%%%%%%%%%%%%%%%%%%%%%%%%%%%%%DIAGRAMME
\put(5,0){\line(1,0){15}}\put(5,5){\line(1,0){20}}\put(0,10){\line(1,0){30}}
\put(0,15){\line(1,0){30}}\put(0,20){\line(1,0){30}}\put(10,25){\line(1,0){10}}
\put(0,10){\line(0,1){10}}\put(5,0){\line(0,1){20}}\put(10,0){\line(0,1){25}}
\put(15,0){\line(0,1){25}}\put(20,0){\line(0,1){25}}\put(25,5){\line(0,1){15}}
\put(30,10){\line(0,1){10}}
\end{picture}}\hspace{2cm}{\setlength{\unitlength}{1mm}
\begin{picture}(35,25)(0,0)
%%%%%%%%%%%%%%%%%%%%%%%%%%%%%%%%%%%%%%%%coloriage
\put(5,0){\color{red}{\rule{5mm}{20mm}}}
%%%%%%%%%%%%%%%%%%%%%%%%%%%%%%%%%%%%%%%%DIAGRAMME
\put(5,0){\line(1,0){15}}\put(5,5){\line(1,0){20}}\put(0,10){\line(1,0){30}}
\put(0,15){\line(1,0){30}}\put(0,20){\line(1,0){30}}\put(10,25){\line(1,0){10}}
\put(0,10){\line(0,1){10}}\put(5,0){\line(0,1){20}}\put(10,0){\line(0,1){25}}
\put(15,0){\line(0,1){25}}\put(20,0){\line(0,1){25}}\put(25,5){\line(0,1){15}}
\put(30,10){\line(0,1){10}}
\end{picture}}
\end{center}
Then we have:
\begin{itemize}
\item $R_{i_1}\prec R_{i_2}\prec R_{i_3}\prec R_{i_4} \prec R_{i_5}$ with $i_1=1$, $i_2=5$, $i_3=4$,
$i_4=2$, $i_5=3$.
\item The column labeled $2$ intersect the rows labeled $2,3,4,5$, thus $a_1=0$,
$a_2=a_3=a_4=a_5=1$.
\end{itemize}
Suppose ${\bf m}=(1,2,1,0,1)$. We then have
\begin{align*}
h_{i_1}=h_1&=r_{1}-a_{1}=2,\\
h_{i_2}=h_5&=r_{5}-m_{1}-a_{5}=1,\\
h_{i_3}=h_4&=r_{4}-(m_{1}+m_{5})-a_{4}=1,\\
h_{i_4}=h_2&=r_{2}-(m_{1}+m_{5}+m_{4})-a_{2}=3,\\
h_{i_5}=h_3&=r_{3}-(m_{1}+m_{5}+m_{4}+m_{2})-a_{3}=1.
\end{align*}
It then follows from Theorem~\ref{thm:distribution} that
$$\sum_{F\in\,\N^c(T,{\bf
m};A)}p^{{\nec}_2(F)}q^{{\sec}_2(F)}=\prod_{j=1}^5{h_j\brack
m_j}_{p,q}={2\brack 1}_{p,q}{3\brack 2}_{p,q}{1\brack
1}_{p,q}{1\brack 0}_{p,q}{1\brack 1}_{p,q}=p^3+2p^2q+2pq^2+q^3.$$

On the other hand, the fillings in $\N^c(T,{\bf m}; A)$ and the
corresponding values of ${\nec}_2$ and ${\sec}_2$ are listed below.

\begin{center}
%%%%%%%%%%%%%%%%%%%%%%%%%%%%%%%%%%%%%%%%%%%%%%%%%%%%diag1
 {\setlength{\unitlength}{1mm}
\begin{picture}(35,25)(0,0)
%%%%%%%%%%%%%%%%%%%%%%%%%%%%%%%%%%%%%%%%coloriage
\put(5,0){\color{red}{\rule{5mm}{20mm}}}
%%%%%%%%%%%%%%%%%%%%%%%%%%%%%%%%%%%%%%%%DIAGRAMME
\put(5,0){\line(1,0){15}}\put(5,5){\line(1,0){20}}\put(0,10){\line(1,0){30}}
\put(0,15){\line(1,0){30}}\put(0,20){\line(1,0){30}}\put(10,25){\line(1,0){10}}
\put(0,10){\line(0,1){10}}\put(5,0){\line(0,1){20}}\put(10,0){\line(0,1){25}}
\put(15,0){\line(0,1){25}}\put(20,0){\line(0,1){25}}\put(25,5){\line(0,1){15}}
\put(30,10){\line(0,1){10}}
%%%%%%%%%%%%%%%%%%%%%%%%%%%%%%%%%%%%%%%%REMPLISSAGE
\put(10,20){\makebox(6,4)[c]{\small\1}}
\put(0,15){\makebox(6,4)[c]{\small\1}}
\put(20,15){\makebox(6,4)[c]{\small \1}}
\put(25,10){\makebox(6,4)[c]{\small \1}}
\put(15,0){\makebox(6,4)[c]{\small \1}}
\put(5,-7){\makebox(15,4)[c]{ ${\nec}_2=0\;,\;{\sec}_2=3$}}
\end{picture}}
\hspace{1cm}
%%%%%%%%%%%%%%%%%%%%%%%%%%%%%%%%%%%%%%%%%%%%%%%%%%%%diag2
  {\setlength{\unitlength}{1mm}
\begin{picture}(35,25)(0,0)
%%%%%%%%%%%%%%%%%%%%%%%%%%%%%%%%%%%%%%%%coloriage
\put(5,0){\color{red}{\rule{5mm}{20mm}}}
%%%%%%%%%%%%%%%%%%%%%%%%%%%%%%%%%%%%%%%%DIAGRAMME
\put(5,0){\line(1,0){15}}\put(5,5){\line(1,0){20}}\put(0,10){\line(1,0){30}}
\put(0,15){\line(1,0){30}}\put(0,20){\line(1,0){30}}\put(10,25){\line(1,0){10}}
\put(0,10){\line(0,1){10}}\put(5,0){\line(0,1){20}}\put(10,0){\line(0,1){25}}
\put(15,0){\line(0,1){25}}\put(20,0){\line(0,1){25}}\put(25,5){\line(0,1){15}}
\put(30,10){\line(0,1){10}}
%%%%%%%%%%%%%%%%%%%%%%%%%%%%%%%%%%%%%%%%REMPLISSAGE
\put(10,20){\makebox(6,4)[c]{\small\1}}
\put(0,15){\makebox(6,4)[c]{\small\1}}
\put(25,15){\makebox(6,4)[c]{\small \1}}
\put(20,10){\makebox(6,4)[c]{\small \1}}
\put(15,0){\makebox(6,4)[c]{\small \1}}
\put(5,-7){\makebox(15,4)[c]{ ${\nec}_2=1\;,\;{\sec}_2=2$}}
\end{picture}}
\hspace{1cm}
%%%%%%%%%%%%%%%%%%%%%%%%%%%%%%%%%%%%%%%%%%%%%%%%%%%%diag3
  {\setlength{\unitlength}{1mm}
\begin{picture}(35,25)(0,0)
%%%%%%%%%%%%%%%%%%%%%%%%%%%%%%%%%%%%%%%%coloriage
\put(5,0){\color{red}{\rule{5mm}{20mm}}}
%%%%%%%%%%%%%%%%%%%%%%%%%%%%%%%%%%%%%%%%DIAGRAMME
\put(5,0){\line(1,0){15}}\put(5,5){\line(1,0){20}}\put(0,10){\line(1,0){30}}
\put(0,15){\line(1,0){30}}\put(0,20){\line(1,0){30}}\put(10,25){\line(1,0){10}}
\put(0,10){\line(0,1){10}}\put(5,0){\line(0,1){20}}\put(10,0){\line(0,1){25}}
\put(15,0){\line(0,1){25}}\put(20,0){\line(0,1){25}}\put(25,5){\line(0,1){15}}
\put(30,10){\line(0,1){10}}
%%%%%%%%%%%%%%%%%%%%%%%%%%%%%%%%%%%%%%%%REMPLISSAGE
\put(10,20){\makebox(6,4)[c]{\small\1}}
\put(20,15){\makebox(6,4)[c]{\small\1}}
\put(25,15){\makebox(6,4)[c]{\small \1}}
\put(0,10){\makebox(6,4)[c]{\small \1}}
\put(15,0){\makebox(6,4)[c]{\small \1}}
\put(5,-7){\makebox(15,4)[c]{ ${\nec}_2=2\;,\;{\sec}_2=1$}}
\end{picture}}
\end{center}
\vspace{1cm}
\begin{center}
%%%%%%%%%%%%%%%%%%%%%%%%%%%%%%%%%%%%%%%%%%%%%%%%%%%%diag4
 {\setlength{\unitlength}{1mm}
\begin{picture}(35,25)(0,0)
%%%%%%%%%%%%%%%%%%%%%%%%%%%%%%%%%%%%%%%%coloriage
\put(5,0){\color{red}{\rule{5mm}{20mm}}}
%%%%%%%%%%%%%%%%%%%%%%%%%%%%%%%%%%%%%%%%DIAGRAMME
\put(5,0){\line(1,0){15}}\put(5,5){\line(1,0){20}}\put(0,10){\line(1,0){30}}
\put(0,15){\line(1,0){30}}\put(0,20){\line(1,0){30}}\put(10,25){\line(1,0){10}}
\put(0,10){\line(0,1){10}}\put(5,0){\line(0,1){20}}\put(10,0){\line(0,1){25}}
\put(15,0){\line(0,1){25}}\put(20,0){\line(0,1){25}}\put(25,5){\line(0,1){15}}
\put(30,10){\line(0,1){10}}
%%%%%%%%%%%%%%%%%%%%%%%%%%%%%%%%%%%%%%%%REMPLISSAGE
\put(15,20){\makebox(6,4)[c]{\small\1}}
\put(0,15){\makebox(6,4)[c]{\small\1}}
\put(20,15){\makebox(6,4)[c]{\small \1}}
\put(25,10){\makebox(6,4)[c]{\small \1}}
\put(10,0){\makebox(6,4)[c]{\small \1}}
\put(5,-7){\makebox(15,4)[c]{ ${\nec}_2=1\;,\;{\sec}_2=2$}}
\end{picture}}
\hspace{1cm}
%%%%%%%%%%%%%%%%%%%%%%%%%%%%%%%%%%%%%%%%%%%%%%%%%%%%diag5
  {\setlength{\unitlength}{1mm}
\begin{picture}(35,25)(0,0)
%%%%%%%%%%%%%%%%%%%%%%%%%%%%%%%%%%%%%%%%coloriage
\put(5,0){\color{red}{\rule{5mm}{20mm}}}
%%%%%%%%%%%%%%%%%%%%%%%%%%%%%%%%%%%%%%%%DIAGRAMME
\put(5,0){\line(1,0){15}}\put(5,5){\line(1,0){20}}\put(0,10){\line(1,0){30}}
\put(0,15){\line(1,0){30}}\put(0,20){\line(1,0){30}}\put(10,25){\line(1,0){10}}
\put(0,10){\line(0,1){10}}\put(5,0){\line(0,1){20}}\put(10,0){\line(0,1){25}}
\put(15,0){\line(0,1){25}}\put(20,0){\line(0,1){25}}\put(25,5){\line(0,1){15}}
\put(30,10){\line(0,1){10}}
%%%%%%%%%%%%%%%%%%%%%%%%%%%%%%%%%%%%%%%%REMPLISSAGE
\put(15,20){\makebox(6,4)[c]{\small\1}}
\put(0,15){\makebox(6,4)[c]{\small\1}}
\put(25,15){\makebox(6,4)[c]{\small \1}}
\put(20,10){\makebox(6,4)[c]{\small \1}}
\put(10,0){\makebox(6,4)[c]{\small \1}}
\put(5,-7){\makebox(15,4)[c]{ ${\nec}_2=2\;,\;{\sec}_2=1$}}
\end{picture}}
\hspace{1cm}
%%%%%%%%%%%%%%%%%%%%%%%%%%%%%%%%%%%%%%%%%%%%%%%%%%%%diag6
  {\setlength{\unitlength}{1mm}
\begin{picture}(35,25)(0,0)
%%%%%%%%%%%%%%%%%%%%%%%%%%%%%%%%%%%%%%%%coloriage
\put(5,0){\color{red}{\rule{5mm}{20mm}}}
%%%%%%%%%%%%%%%%%%%%%%%%%%%%%%%%%%%%%%%%DIAGRAMME
\put(5,0){\line(1,0){15}}\put(5,5){\line(1,0){20}}\put(0,10){\line(1,0){30}}
\put(0,15){\line(1,0){30}}\put(0,20){\line(1,0){30}}\put(10,25){\line(1,0){10}}
\put(0,10){\line(0,1){10}}\put(5,0){\line(0,1){20}}\put(10,0){\line(0,1){25}}
\put(15,0){\line(0,1){25}}\put(20,0){\line(0,1){25}}\put(25,5){\line(0,1){15}}
\put(30,10){\line(0,1){10}}
%%%%%%%%%%%%%%%%%%%%%%%%%%%%%%%%%%%%%%%%REMPLISSAGE
\put(15,20){\makebox(6,4)[c]{\small\1}}
\put(20,15){\makebox(6,4)[c]{\small\1}}
\put(25,15){\makebox(6,4)[c]{\small \1}}
\put(0,10){\makebox(6,4)[c]{\small \1}}
\put(10,0){\makebox(6,4)[c]{\small \1}}
\put(5,-7){\makebox(15,4)[c]{ ${\nec}_2=3\;,\;{\sec}_2=0$}}
\end{picture}}
\end{center}
\vspace{1cm}

Summing up we get $\sum_{F\in\N^c(T,{\bf m};
A)}p^{{\nec}_{2}(F)}q^{{\sec}_{2}(F)}=p^3+2p^2q+2pq^2+q^3$, as
desired.\\

\section{Symmetry of $2$-crossings and $2$-nestings
in linked partitions,\\ set partitions and matchings}

 In this section we show briefly how results on the enumeration
 of $2$-crossings and $2$-nestings can be recovered from the results
 obtained in this paper. Let $G$ be a simple
graph on $[n]$. The multiset of lefthand (resp., righthand)
endpoints of the arcs of $G$ will be denoted by $\le(G)$ (resp.,
$\ri(G)$). For instance, if $G$ is the graph drawn in
Figure~\ref{fig:part-filling}, we have $\le(G) = \{1,2,2,3,
5,6,6,9\}$ and $\ri(G) = \{3,4, 6, 7, 9,9,10,11\}$. For $S$ and $T$
two multisubsets of $[n]$, we will denote  by $\G_n(S,T)$ the set of
simple graphs $G$ on $[n]$ satisfying $\le(G)=S$ and $\ri(G)=T$.

Suppose $F$ is the $01$-filling of the triangular shape $\Delta_n$
which corresponds to the graph~$G$. For convenience, we joined an
empty column at the right and an empty row at the top of $\Delta_n$,
and columns are labeled from left to right and rows from top to
bottom by~$\{1,2,\ldots,n\}$. It is obvious that the number of $1$'s
in the column (resp., row) labeled~$i$  is equal to the multiplicity
of $i$ in $\le(G)$ (resp., $\ri(G)$). See
Figure~\ref{fig:part-filling} for an illustration. Taking the moon
polyomino $T:=\Delta_n$ in Theorem~\ref{thm:distribution}, we obtain
the following results.

 Let $(S,T)$ be a pair of multisubsets of $[n]$ and denote by
 $m_i$ the multiplicity of $i$ in $T$ and by $m'_i$ the multiplicity of $i$ in $S$.
Also for any $i\in T$ set $h_i=|\{j\in S\,|\,j<i\}|-|\{j\in
T\,|\,j<i\}|$ and for any $i\in S$ set $h'_i=|\{j\in
T\,|\,j>i\}|-|\{j\in S\,|\,j>i\}|$.

\begin{cor}\label{cor:symgraphes}
Let $(S,T)$ be a pair of multisubsets of $[n]$.
\begin{enumerate}
\item
If all elements of $S$ have multiplicity $1$, then
\begin{equation}\label{eq:distgraphe1}
\sum_{G\in\,\G_n(S,T)}p^{{\nec}_2(F)}q^{{\sec}_2(F)}
=\sum_{G\in\,\G_n(S,T)}p^{{\sec}_2(F)}q^{{\nec}_2(F)}=\prod_{i\in
T}{h_i\brack m_i}_{p,q}.
\end{equation}
\item
If all elements of $T$ have multiplicity $1$, then
\begin{equation}\label{eq:distgraphe2}
\sum_{G\in\,\G_n(S,T)}p^{{\nec}_2(F)}q^{{\sec}_2(F)}
=\sum_{G\in\,\G_n(S,T)}p^{{\sec}_2(F)}q^{{\nec}_2(F)}=\prod_{i\in
S}{h'_i\brack m'_i}_{p,q}.
\end{equation}
\end{enumerate}
In particular,  the joint statistic $({\cro}_2,{\ne}_2)$ is
symmetrically distributed over $\G_n(S,T)$ if
\begin{itemize}
\item either all elements of $S$ have multiplicity $1$,
\item either all elements of $T$ have multiplicity $1$.
\end{itemize}
\end{cor}

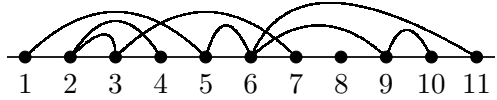
\begin{figure}[h!]
\begin{center}
{\setlength{\unitlength}{1.2mm}
\begin{picture}(50,10)(0,-3)
\put(-2,0){\line(1,0){54}}
\put(0,0){\circle*{1,3}}\put(0,0){\makebox(0,-6)[c]{\small $1$}}
\put(5,0){\circle*{1,3}}\put(5,0){\makebox(0,-6)[c]{\small $2$}}
\put(10,0){\circle*{1,3}}\put(10,0){\makebox(0,-6)[c]{\small $3$}}
\put(15,0){\circle*{1,3}}\put(15,0){\makebox(0,-6)[c]{\small $4$}}
\put(20,0){\circle*{1,3}}\put(20,0){\makebox(0,-6)[c]{\small $5$}}
\put(25,0){\circle*{1,3}}\put(25,0){\makebox(0,-6)[c]{\small $6$}}
\put(30,0){\circle*{1,3}}\put(30,0){\makebox(0,-6)[c]{\small $7$}}
\put(35,0){\circle*{1,3}}\put(35,0){\makebox(0,-6)[c]{\small $8$}}
\put(40,0){\circle*{1,3}}\put(40,0){\makebox(0,-6)[c]{\small $9$}}
\put(45,0){\circle*{1,3}}\put(45,0){\makebox(0,-6)[c]{\small $10$}}
\put(50,0){\circle*{1,3}}\put(50,0){\makebox(0,-6)[c]{\small $11$}}
\qbezier(0,0)(10,10)(20,0)
\qbezier(5,0)(10,5)(10,0)\qbezier(20,0)(22,7)(25,0)
\qbezier(10,0)(20,10)(30,0)\qbezier(25,0)(30,12)(50,0)
\qbezier(25,0)(32.5,7)(40,0)\qbezier(40,0)(42,6)(45,0)
\qbezier(5,0)(10,8)(15,0)
\end{picture}}
\end{center}
\caption{Linear representation of
$\pi=\{1,5\}\{2,3,4\}\{3,7\}\{5,6\}\{6,9,11\}\{8\}\{9,10\}$}\label{fig:linkedpartition}
\end{figure}

Note that \eqref{eq:distgraphe2} is equivalent to a result of Chen
et al.~\cite[Theorem 3.5]{Chen1} on the enumeration of $2$-crossings
and $2$-nestings in linked set partitions. Let $E$ and $F$ be two
finite sets of nonnegative integers. We say that $E$ and $F$ are
\emph{nearly disjoint} if for every $i\in E\cap F$, one of the
following holds:
\begin{itemize}
  \item[(a)] $i = \min(E)$, $|E| > 1$ and $i\neq \min(F)$, or
  \item[(b)] $i = \min(F)$, $|F| > 1$ and $i\neq \min(E)$.
\end{itemize}
 A \emph{linked partition} (see \cite{Chen1}) of
 $[n]$ is a collection of non-empty and
pairwise nearly disjoint subsets whose union is $[n]$.
 The set of all linked partitions
of $[n]$ will be denoted by $\LP_n$. The linear representation
$G_{\pi}$ of a linked partition $\pi\in\LP_n$ is the graph on $[n]$
where $i$ and $j$ are connected by an arc if and only if $j$ lies in
a block $B$ with $i= {\Min}(B)$. An illustration is given in
Figure~\ref{fig:linkedpartition}.
Clearly, the map $\pi\mapsto G_{\pi}$ establishes a bijection
between linked set partitions and simple graphs $G$ such that all
elements of $\ri(G)$ have multiplicity one. For $S,T \subseteq[n]$
two multisubsets of $[n]$, denote by $\LP_n(S,T)$ the set
$\{\pi\in\LP_n\,:\, \le(G_{\pi})= S\,,\,\ri(G_{\pi}) = T\}$. Then
\eqref{eq:distgraphe2} can be rewritten
\begin{equation}\label{eq:sym-linked}
\sum_{\pi\in\,\LP_n(S,T)}p^{{\cro}_2(G_{\pi})}q^{{\ne}_2(G_{\pi})}
=\sum_{\pi\in\,\LP_n(S,T)}p^{{\ne}_2(G_{\pi})}q^{{\cro}_2(G_{\pi})}
=\prod_{i\in S}{h'_i\brack m'_i}_{p,q},
\end{equation}
where $h'_i$ and $m'_i$ are defined as in
Corollary~\eqref{cor:symgraphes}, which is equivalent to a result of
Chen et al.~\cite[Theorem 3.5]{Chen1}.

Now consider the map $\pi\mapsto St_{\pi}$ which associates to each
set partition its standard representation (see
Figure~\ref{fig:partition-matching}). Clearly, this map establishes
a bijection between set partitions and simple graphs $G$ such that
all elements of $\ri(G)$ and $\le(G)$ have multiplicity one.
Applying Corollary~\ref{cor:symgraphes} with $m_i=1$ (or
\eqref{eq:sym-linked} with $m'_i=1$) we
recover~\eqref{eq:sym-partition} and the following identity which is
implicit in~\cite[Section~4]{KaZe}
\begin{equation}
\sum_{\pi\in\,\P_n(S,T)}p^{{\cro}_2(\pi)}q^{{\ne}_2(\pi)}
=\sum_{\pi\in\,\P_n(S,T)}p^{{\ne}_2(\pi)}q^{{\cro}_2(\pi)}
=\prod_{i\in\O}[h_i]_{p,q}=\prod_{i\in T}[h'_i]_{p,q},
\end{equation}
where $h_i$ and $h'_i$ are defined as in
Corollary~\eqref{cor:symgraphes}.

\section{Proof of Theorem~\ref{thm:distribution}}
As explained in Section~2, it suffices to prove the first part of
Theorem~\ref{thm:distribution} that is the
identity~\eqref{eq:symfilling}. Throughout this section, $T$ is a
moon polyomino with $s$ rows and length-row sequence
$(r_1,r_2,\ldots,r_s)$, and $\textbf{m}=(m_1,\ldots,m_s)$ is a
$s$-uple of positive integers.

\subsection{Preliminaries}

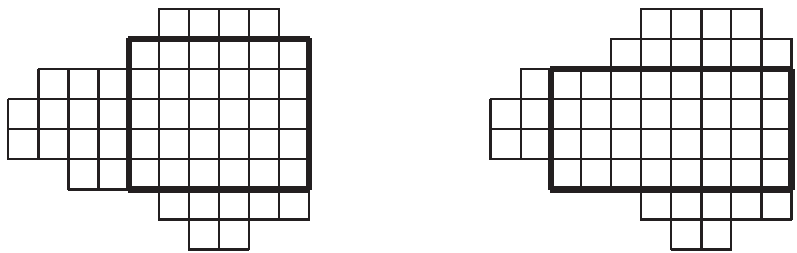
\begin{figure}[h!]
\begin{center}
 {\setlength{\unitlength}{0.8mm}
\begin{picture}(50,50)(0,0)
%%%%%%%%%%%%%%%%%%%%%%%%%%%%%%%%%%%%%%%%DIAGRAMME
\put(30,0){\line(1,0){10}} \put(25,5){\line(1,0){25}}
\put(10,10){\line(1,0){40}} \put(0,15){\line(1,0){50}}
\put(0,20){\line(1,0){50}} \put(0,25){\line(1,0){50}}
\put(5,30){\line(1,0){45}} \put(20,35){\line(1,0){30}}
\put(15,30){\line(1,0){10}} \put(25,40){\line(1,0){20}}
\put(0,15){\line(0,1){10}}\put(5,15){\line(0,1){15}}
\put(10,10){\line(0,1){20}}\put(15,10){\line(0,1){20}}\put(20,10){\line(0,1){25}}
\put(25,5){\line(0,1){35}}\put(30,0){\line(0,1){40}}\put(35,0){\line(0,1){40}}\put(40,0){\line(0,1){40}}
\put(45,5){\line(0,1){35}} \put(50,5){\line(0,1){30}}
%%%%%%%%%%%%%%%%%%%%%%%%%%%%%%%%%%%%%%%%
\linethickness{0.7mm}
\put(20,35){\line(1,0){30}}\put(20,10){\line(1,0){30}}
\put(20,10){\line(0,1){25}}\put(50,10){\line(0,1){25}}
\end{picture}}
\hspace{2cm} {\setlength{\unitlength}{0.8mm}
\begin{picture}(50,50)(0,0)
%%%%%%%%%%%%%%%%%%%%%%%%%%%%%%%%%%%%%%%%DIAGRAMME
\put(30,0){\line(1,0){10}} \put(25,5){\line(1,0){25}}
\put(10,10){\line(1,0){40}} \put(0,15){\line(1,0){50}}
\put(0,20){\line(1,0){50}} \put(0,25){\line(1,0){50}}
\put(5,30){\line(1,0){45}} \put(20,35){\line(1,0){30}}
\put(15,30){\line(1,0){10}} \put(25,40){\line(1,0){20}}
\put(0,15){\line(0,1){10}}\put(5,15){\line(0,1){15}}
\put(10,10){\line(0,1){20}}\put(15,10){\line(0,1){20}}\put(20,10){\line(0,1){25}}
\put(25,5){\line(0,1){35}}\put(30,0){\line(0,1){40}}\put(35,0){\line(0,1){40}}\put(40,0){\line(0,1){40}}
\put(45,5){\line(0,1){35}} \put(50,5){\line(0,1){30}}
%%%%%%%%%%%%%%%%%%%%%%%%%%%%%%%%%%%%%%%%
\linethickness{0.7mm}
\put(10,10){\line(1,0){40}}\put(10,30){\line(1,0){40}}
\put(10,10){\line(0,1){20}}\put(50,10){\line(0,1){20}}
\end{picture}}
\end{center}
\caption{\emph{left}: the $2$-th rectangle, \emph{right}: the $6$-th
rectangle.}\label{fig:rectanglepolyomino}
\end{figure}

Let $i$, $1\leq i\leq s$, be an integer. The \emph{$i$-th rectangle
of $T$}, is the greatest rectangle contained in $T$ whose top
(resp., bottom) row is $R_i$ if $R_i\in Up(T)$ (resp., $R_i\in
Low(T)$). An illustration is given in
Figure~\ref{fig:rectanglepolyomino}.

Let $F$ be a $01$-filling of $T$ in $\N^c(T,{\bf m})$. The
\emph{coloring of $F$} is the colored filling obtained from $F$ by:
\begin{itemize}
\item coloring the cells of the empty columns,
\item for $i=1,\ldots,s$, coloring the cells which are both contained in
the $i$-th rectangle and
\begin{itemize}
\item if $R_i\in Up(T)$, below
a $1$ of $R_i$.
\item if $R_i\in Low(T)$, above a
$1$ of $R_i$.
\end{itemize}
\end{itemize}
An illustration is given in Figure~\ref{fig:coloring of a moon
polyomino}. Throughout this paper, we identify a filling with its
coloring. For instance, "the cell $c$ of the filling $F$ is
uncolored" means that "the cell $c$ is uncolored in the coloring of
$F$",\ldots

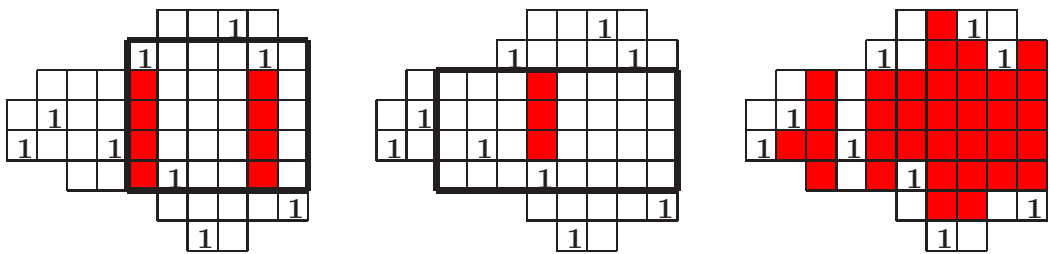
\begin{figure}[h!]
\begin{center}
 {\setlength{\unitlength}{0.8mm}
\begin{picture}(50,50)(0,0)
%%%%%%%%%%%%%%%%%%%%%%%%%%%%%%%%%%%%%%%%coloriage
\put(20,10){\color{red}{\rule{4mm}{16mm}}}\put(40,10){\color{red}{\rule{4mm}{16mm}}}
%%%%%%%%%%%%%%%%%%%%%%%%%%%%%%%%%%%%%%%%DIAGRAMME
\put(30,0){\line(1,0){10}} \put(25,5){\line(1,0){25}}
\put(10,10){\line(1,0){40}} \put(0,15){\line(1,0){50}}
\put(0,20){\line(1,0){50}} \put(0,25){\line(1,0){50}}
\put(5,30){\line(1,0){45}} \put(20,35){\line(1,0){30}}
\put(15,30){\line(1,0){10}} \put(25,40){\line(1,0){20}}
\put(0,15){\line(0,1){10}}\put(5,15){\line(0,1){15}}
\put(10,10){\line(0,1){20}}\put(15,10){\line(0,1){20}}\put(20,10){\line(0,1){25}}
\put(25,5){\line(0,1){35}}\put(30,0){\line(0,1){40}}\put(35,0){\line(0,1){40}}\put(40,0){\line(0,1){40}}
\put(45,5){\line(0,1){35}} \put(50,5){\line(0,1){30}}
%%%%%%%%%%%%%%%%%%%%%%%%%%%%%%%%%%%%%%%%REMPLISSAGE
\put(30,0){\makebox(6,4)[c]{\small \1}}
\put(45,5){\makebox(6,4)[c]{\small \1}}
\put(25,10){\makebox(6,4)[c]{\small \1}}
\put(15,15){\makebox(6,4)[c]{\small\1}}
\put(0,15){\makebox(6,4)[c]{\small \1}}
\put(5,20){\makebox(6,4)[c]{\small\1}}
\put(20,30){\makebox(6,4)[c]{\small\1}}
\put(40,30){\makebox(6,4)[c]{\small \1}}
\put(35,35){\makebox(6,4)[c]{\small \1}}
%%%%%%%%%%%%%%%%%%%%%%%%%%%%%%%%%%%%%%%%
\linethickness{0.7mm}
\put(20,35){\line(1,0){30}}\put(20,10){\line(1,0){30}}
\put(20,10){\line(0,1){25}}\put(50,10){\line(0,1){25}}
\end{picture}}
\hspace{0.5cm} {\setlength{\unitlength}{0.8mm}
\begin{picture}(50,50)(0,0)
%%%%%%%%%%%%%%%%%%%%%%%%%%%%%%%%%%%%%%%%coloriage
\put(25,15){\color{red}{\rule{4mm}{12mm}}}
%%%%%%%%%%%%%%%%%%%%%%%%%%%%%%%%%%%%%%%%DIAGRAMME
\put(30,0){\line(1,0){10}} \put(25,5){\line(1,0){25}}
\put(10,10){\line(1,0){40}} \put(0,15){\line(1,0){50}}
\put(0,20){\line(1,0){50}} \put(0,25){\line(1,0){50}}
\put(5,30){\line(1,0){45}} \put(20,35){\line(1,0){30}}
\put(15,30){\line(1,0){10}} \put(25,40){\line(1,0){20}}
\put(0,15){\line(0,1){10}}\put(5,15){\line(0,1){15}}
\put(10,10){\line(0,1){20}}\put(15,10){\line(0,1){20}}\put(20,10){\line(0,1){25}}
\put(25,5){\line(0,1){35}}\put(30,0){\line(0,1){40}}\put(35,0){\line(0,1){40}}\put(40,0){\line(0,1){40}}
\put(45,5){\line(0,1){35}} \put(50,5){\line(0,1){30}}
%%%%%%%%%%%%%%%%%%%%%%%%%%%%%%%%%%%%%%%%REMPLISSAGE
\put(30,0){\makebox(6,4)[c]{\small \1}}
\put(45,5){\makebox(6,4)[c]{\small \1}}
\put(25,10){\makebox(6,4)[c]{\small \1}}
\put(15,15){\makebox(6,4)[c]{\small\1}}
\put(0,15){\makebox(6,4)[c]{\small \1}}
\put(5,20){\makebox(6,4)[c]{\small\1}}
\put(20,30){\makebox(6,4)[c]{\small\1}}
\put(40,30){\makebox(6,4)[c]{\small \1}}
\put(35,35){\makebox(6,4)[c]{\small \1}}
%%%%%%%%%%%%%%%%%%%%%%%%%%%%%%%%%%%%%%%%
\linethickness{0.7mm}
\put(10,10){\line(1,0){40}}\put(10,30){\line(1,0){40}}
\put(10,10){\line(0,1){20}}\put(50,10){\line(0,1){20}}
\end{picture}}
\hspace{0.5cm}
 {\setlength{\unitlength}{0.8mm}
\begin{picture}(50,50)(0,0)
%%%%%%%%%%%%%%%%%%%%%%%%%%%%%%%%%%%%%%%%coloriage
\put(30,5){\color{red}{\rule{4mm}{28mm}}}\put(35,5){\color{red}{\rule{4mm}{24mm}}}
\put(45,10){\color{red}{\rule{4mm}{20mm}}}\put(25,15){\color{red}{\rule{4mm}{12mm}}}
\put(5,15){\color{red}{\rule{4mm}{4mm}}}\put(10,10){\color{red}{\rule{4mm}{16mm}}}
\put(20,10){\color{red}{\rule{4mm}{16mm}}}\put(40,10){\color{red}{\rule{4mm}{16mm}}}
%%%%%%%%%%%%%%%%%%%%%%%%%%%%%%%%%%%%%%%%DIAGRAMME
\put(30,0){\line(1,0){10}} \put(25,5){\line(1,0){25}}
\put(10,10){\line(1,0){40}} \put(0,15){\line(1,0){50}}
\put(0,20){\line(1,0){50}} \put(0,25){\line(1,0){50}}
\put(5,30){\line(1,0){45}} \put(20,35){\line(1,0){30}}
\put(15,30){\line(1,0){10}} \put(25,40){\line(1,0){20}}
\put(0,15){\line(0,1){10}}\put(5,15){\line(0,1){15}}
\put(10,10){\line(0,1){20}}\put(15,10){\line(0,1){20}}\put(20,10){\line(0,1){25}}
\put(25,5){\line(0,1){35}}\put(30,0){\line(0,1){40}}\put(35,0){\line(0,1){40}}\put(40,0){\line(0,1){40}}
\put(45,5){\line(0,1){35}} \put(50,5){\line(0,1){30}}
%%%%%%%%%%%%%%%%%%%%%%%%%%%%%%%%%%%%%%%%REMPLISSAGE
\put(30,0){\makebox(6,4)[c]{\small \1}}
\put(45,5){\makebox(6,4)[c]{\small \1}}
\put(25,10){\makebox(6,4)[c]{\small \1}}
\put(15,15){\makebox(6,4)[c]{\small\1}}
\put(0,15){\makebox(6,4)[c]{\small \1}}
\put(5,20){\makebox(6,4)[c]{\small\1}}
\put(20,30){\makebox(6,4)[c]{\small\1}}
\put(40,30){\makebox(6,4)[c]{\small \1}}
\put(35,35){\makebox(6,4)[c]{\small \1}}
\end{picture}}
\end{center}
\caption{\emph{left}: coloring induced by $R_2$, \emph{center}:
coloring induced by $R_6$, \emph{right}: full
coloring.}\label{fig:coloring of a moon polyomino}
\end{figure}

The interest of coloring a $01$-filling is in the following result.
Let $c$ be a cell of $F$. If $c$ contains a $1$ denote by
${\luc}(c;F)$ (resp., ${\ruc}(c;F)$) the numbers of uncolored empty
cells which are both to the left (resp., right) and in the same row
than the cell $c$ in $F$. If $c$ is empty, set
${\luc}(c;F)={\ruc}(c;F)=0$.

\begin{prop}\label{prop:colored filling}
Let $F\in\N^c(T)$ and $c$ be a cell of $R_i$ containing a 1. Then
${\luc}(c;F)$ (resp., ${\ruc}(c;F)$) is equal to
\begin{itemize}
\item if $R_i\in Up(T)$: the number of ascents (resp., descents) contained in the $i$-th rectangle of $F$ whose
North-east (resp., North-west) $1$ is in $c$,
\item if $R_i\in Low(T)$: the number of descents (resp. ascents)  contained in the $i$-th rectangle of $F$ whose South-east
(resp., South-west) $1$ is in $c$.
\end{itemize}
\end{prop}

\emph{Sketch of the proof of Proposition~\ref{prop:colored
filling}.} Suppose $R_i\in Up(T)$ and let $c$ be a cell of~$R_i$
containing a~$1$.

Let $c'$ be an empty uncolored cell in $R_i$ to the left (resp.,
right) of $c$. Suppose $c'$ belong to the column $C_k$. By the
definition of the coloring of polyominoes, the column $C_k$ contains
a 1 (otherwise all the cells of $C_k$, in particular $c'$, would be
colored). Moreover, the cell $c''$ of $C_k$ containing a $1$ must
belong to a row $R_j$ with $R_i\prec R_j$ (otherwise all the cells
of $C_k$ in the $i$-th rectangle of $T$, in particular $c'$, would
be colored), and thus $c''$ belong to the $i$-th rectangle. Since
$R_i\in Up(T)$, the row $R_i$ is the top row of the $i$-th rectangle
of $T$, and thus the cell $c''$ is to the South-west (resp.,
South-east) of the cell $c$. Finally, the sequence $c''c$ (resp.,
$cc''$ ) is an ascent (resp., a descent)  contained in the $i$-th
rectangle of~$F$.

Reversely, let $c''$ be a cell of $F$ such that the pair $c''c$
(resp., $cc''$) is an ascent (resp., descent) of $F$ contained in
the $i$-th rectangle of $F$. Suppose $c''$ belong to $C_k$ and let
$c'$ be the cell of $F$ at the intersection of the column $C_k$ and
the row $R_i$. Clearly, $C_k$ is empty (there is at most one 1 in
each column). It remains to show that the cell $c'$ is uncolored.
This follows from the fact that $c''$ belong to a row $R_j$ with
$R_i\prec R_j$ (since $c''$ is in the $i$-th rectangle). We thus
have proved the first part of Proposition~\ref{prop:colored
filling}.

The second part can be proved by a similar reasoning. Therefore, the
details are left  to the reader.

 \qed

Note that Proposition~\ref{prop:colored filling} lead to the
following decompositions of ${\nec}_2$ and ${\sec}_2$:
\begin{align}
{\nec}_2(F)=\sum_{c \,\in \,Up(F)} {\luc}(c;F)+\sum_{c\, \in \,Low(F)}{\ruc}(c;F),\label{eq:cro-sum}\\
{\sec}_2(F)=\sum_{c \,\in \,Up(F)}{\ruc}(c;F)+\sum_{c
\,\in\,Low(F)}{\luc}(c;F).\label{eq:ne-sum}
\end{align}

\subsection{A correspondence between $01$-fillings and sequence of compositions}

If $n$ and $k$ are positive integers, we will denote by $\C_k(n)$
the set of compositions of $n$ into $k$ positive parts. Recall that
a element in $\C_k(n)$ is just a $k$-uple $(b_1,b_2,\ldots,b_k)$ of
positive integers such that $b_1+b_2+\cdots+b_k=n$. The proof of
Theorem~\ref{thm:distribution} is based on a bijection
$$\Psi:\N^c(T,{\bf m};A)\to
\C_{m_{1}+1}(h_{1}-m_1)\times\C_{m_{2}+1}(h_{2}-m_2)\times
\cdots\times \C_{m_{s}+1}(h_{s}-m_s)$$ which keeps track of the
statistics ${\nec}_2$ and ${\sec}_2$.

\subsubsection*{Algorithm for $\Psi$}
\begin{figure}[h!]
\begin{center}
 {\setlength{\unitlength}{0.8mm}
\begin{picture}(50,50)(0,-5)
\put(5,35){\makebox(6,4)[c]{$F$}}
%%%%%%%%%%%%%%%%%%%%%%%%%%%%%%%%%%%%%%%%coloriage
\put(30,5){\color{red}{\rule{4mm}{28mm}}}\put(35,5){\color{red}{\rule{4mm}{24mm}}}
\put(45,10){\color{red}{\rule{4mm}{20mm}}}\put(25,15){\color{red}{\rule{4mm}{12mm}}}
\put(5,15){\color{red}{\rule{4mm}{4mm}}}\put(10,10){\color{red}{\rule{4mm}{16mm}}}
\put(20,10){\color{red}{\rule{4mm}{16mm}}}\put(40,10){\color{red}{\rule{4mm}{16mm}}}
%%%%%%%%%%%%%%%%%%%%%%%%%%%%%%%%%%%%%%%%DIAGRAMME
\put(30,0){\line(1,0){10}} \put(25,5){\line(1,0){25}}
\put(10,10){\line(1,0){40}} \put(0,15){\line(1,0){50}}
\put(0,20){\line(1,0){50}} \put(0,25){\line(1,0){50}}
\put(5,30){\line(1,0){45}} \put(20,35){\line(1,0){30}}
\put(15,30){\line(1,0){10}} \put(25,40){\line(1,0){20}}
\put(0,15){\line(0,1){10}}\put(5,15){\line(0,1){15}}
\put(10,10){\line(0,1){20}}\put(15,10){\line(0,1){20}}\put(20,10){\line(0,1){25}}
\put(25,5){\line(0,1){35}}\put(30,0){\line(0,1){40}}\put(35,0){\line(0,1){40}}\put(40,0){\line(0,1){40}}
\put(45,5){\line(0,1){35}} \put(50,5){\line(0,1){30}}
%%%%%%%%%%%%%%%%%%%%%%%%%%%%%%%%%%%%%%%%REMPLISSAGE
\put(30,0){\makebox(6,4)[c]{\small \1}}
\put(45,5){\makebox(6,4)[c]{\small \1}}
\put(25,10){\makebox(6,4)[c]{\small \1}}
\put(15,15){\makebox(6,4)[c]{\small\1}}
\put(0,15){\makebox(6,4)[c]{\small \1}}
\put(5,20){\makebox(6,4)[c]{\small\1}}
\put(20,30){\makebox(6,4)[c]{\small\1}}
\put(40,30){\makebox(6,4)[c]{\small \1}}
\put(35,35){\makebox(6,4)[c]{\small \1}}
\end{picture}}
\hspace{2cm}
 {\setlength{\unitlength}{0.8mm}
\begin{picture}(20,50)(0,-5)
%%%%%%%%%%%%%%%%%%%%%%%%%%%%%%%%%%%%%%%%COMPOSITIONS
\put(0,35){\makebox(6,4)[l]{$\c^{(1)}=(1,1)$}}
\put(0,30){\makebox(6,4)[l]{$\c^{(2)}=(0,1,0)$}}
\put(0,25){\makebox(6,4)[l]{$\c^{(3)}=(0)$}}
\put(0,20){\makebox(6,4)[l]{$\c^{(4)}=(1,1)$}}
\put(0,15){\makebox(6,4)[l]{$\c^{(5)}=(0,0,0)$}}
\put(0,10){\makebox(6,4)[l]{$\c^{(6)}=(1,0)$}}
\put(0,5){\makebox(6,4)[l]{$\c^{(7)}=(2,0)$}}
\put(0,0){\makebox(6,4)[l]{$\c^{(8)}=(0,1)$}}
\put(-20,15){\makebox(6,4)[l]{$\longrightarrow$}}
\put(-18,19){\makebox(6,4)[c]{$\Psi$}}
\end{picture}}
\caption{The mapping $\Psi$}\label{fig:Psi}
\end{center}
\end{figure}
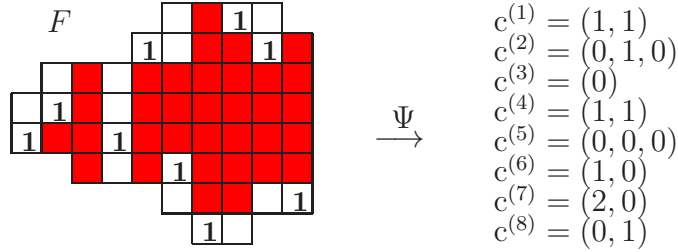

For $F\in\N^c(T,{\bf m};A)$ associate the sequence of compositions
$\Psi(F):=(\c^{(1)},\c^{(2)},\ldots,\c^{(s)})$, where for
$i=1,\ldots,s$, the composition $\c^{(i)}$ is defined by
\begin{itemize}
\item $\c^{(i)}=(0)$ if $m_i=0$, otherwise
\item $\c^{(i)}=(\c^{(i)}_1,\ldots,\c^{(i)}_{m_i+1})$ where $\c^{(i)}_1$ (resp., $\c^{(i)}_j$ for $2\leq j\leq m_i$,
$\c^{(i)}_{m_i+1}$) is the number of uncolored cells in $R_i$ to the
left of the first~$1$ (resp., between the $j$-th 1 and the
$(j+1)$-th~$1$, to the right of the last~$1$) of $R_i$ in the
coloring of~$F$.
\end{itemize}
An illustration is given in Figure~\ref{fig:Psi}.

In order to show that $\Psi$ is bijective, we describe its reverse.
Let ${\textrm{\textbf{c}}}=(\c^{(1)},\c^{(2)},\ldots,\c^{(s)})$ in
$\C_{m_{1}+1}(h_{1}-m_1)\times\C_{m_{2}+1}(h_{2}-m_2)\times
\cdots\times \C_{m_{s}+1}(h_{s}-m_s)$. Then define the $01$-filling
$\Upsilon({\textrm{\textbf{c}}})$ of $T$ by the following process.

 (1) Color the columns indexed by the set $A$ of the polyomino $T$. Denote by $F_0$ the
 colored polyomino obtained.

(2) Construct a sequence of colored fillings $(F_j)_{j=1...s}$ of
$T$ as follows. Suppose $R_{i_1}\prec R_{i_2} \prec \cdots \prec
R_{i_s}$. Then for $j$ from 1 to $s$, the (colored) filling $F_{j}$
is obtained from $F_{j-1}$~as follows:

\begin{itemize}
\item if $m_{i_j}=0$, do nothing,
\item else, insert $m_{i_j}$ 1's in the $i_j\,$-th row of $F_{j-1}$ in such a
way that the number of uncolored cells strictly
\begin{itemize}
\item to the left of the first~$1$ is $\c^{(i_j)}_1$,
\item between the $u$-th 1 and the $(u+1)$-th 1,
$1\leq u\leq m_{i_j}-1$, is $\c^{(i_j)}_{u+1}$,
\item to the right of the last~$1$ is $\c^{(i_j)}_{m_{i_j}+1}$.
\end{itemize}
Next, color the cells which are both below (resp., above) the new
$1$'s inserted in $R_{i_j}$ and contained in the $i_j$-th rectangle
if $R_{i_j}\in Up(T)$ (resp., $R_{i_j}\in Low(T)$).
\end{itemize}

(3) Set $\Upsilon({\textrm{\textbf{c}}})=F_s$.\\

For a better understanding, we give an example. Suppose $T$ is the
moon polyomino given below, $A=\{2\}$ and ${\bf m}=(1,2,1,0,1)$.
Note that $R_1\prec R_5\prec R_4\prec R_2\prec R_3$.

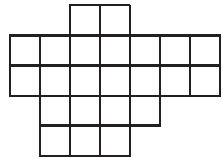
\begin{center}
{\setlength{\unitlength}{0.8mm}
\begin{picture}(35,25)(0,0)
\put(5,0){\line(1,0){15}}\put(5,5){\line(1,0){20}}\put(0,10){\line(1,0){35}}
\put(0,15){\line(1,0){35}}\put(0,20){\line(1,0){35}}\put(10,25){\line(1,0){10}}
\put(0,10){\line(0,1){10}}\put(5,0){\line(0,1){20}}\put(10,0){\line(0,1){25}}
\put(15,0){\line(0,1){25}}\put(20,0){\line(0,1){25}}\put(25,5){\line(0,1){15}}
\put(30,10){\line(0,1){10}}\put(35,10){\line(0,1){10}}
\end{picture}}
\end{center}

\begin{figure}[h!]\label{fig:compo-fillings}
\begin{center}
%%%%%%%%%%%%%%%%%%%%%%%%%%%%%%%%%%%%%%%%%%%%%%%%%%%%%%%%%%%%%%%%%%%%%%%%%%%%%%%%%%%%%%%%%%%%%%%%STEP0
 {\setlength{\unitlength}{1mm}
\begin{picture}(35,25)(0,0)
\put(-2,3){\makebox(6,4)[c]{$F_0$}}
%%%%%%%%%%%%%%%%%%%%%%%%%%%%%%%%%%%%%%%%coloriage
\put(5,0){\color{red}{\rule{5mm}{20mm}}}
%%%%%%%%%%%%%%%%%%%%%%%%%%%%%%%%%%%%%%%%DIAGRAMME
\put(5,0){\line(1,0){15}}\put(5,5){\line(1,0){20}}\put(0,10){\line(1,0){35}}
\put(0,15){\line(1,0){35}}\put(0,20){\line(1,0){35}}\put(10,25){\line(1,0){10}}
\put(0,10){\line(0,1){10}}\put(5,0){\line(0,1){20}}\put(10,0){\line(0,1){25}}
\put(15,0){\line(0,1){25}}\put(20,0){\line(0,1){25}}\put(25,5){\line(0,1){15}}
\put(30,10){\line(0,1){10}}\put(35,10){\line(0,1){10}}
\end{picture}}
\hspace{1cm}
%%%%%%%%%%%%%%%%%%%%%%%%%%%%%%%%%%%%%%%%%%%%%%%%%%%%%%%%%%%%%%%%%%%%%%%%%%%%%%%%%%%%%%%%%%%%%%%%STEP1
 {\setlength{\unitlength}{1mm}
\begin{picture}(35,25)(0,0)
\put(-2,3){\makebox(6,4)[c]{$F_1$}}
%%%%%%%%%%%%%%%%%%%%%%%%%%%%%%%%%%%%%%%%coloriage
\put(5,0){\color{red}{\rule{5mm}{20mm}}}
\put(15,0){\color{red}{\rule{5mm}{20mm}}}
%%%%%%%%%%%%%%%%%%%%%%%%%%%%%%%%%%%%%%%%DIAGRAMME
\put(5,0){\line(1,0){15}}\put(5,5){\line(1,0){20}}\put(0,10){\line(1,0){35}}
\put(0,15){\line(1,0){35}}\put(0,20){\line(1,0){35}}\put(10,25){\line(1,0){10}}
\put(0,10){\line(0,1){10}}\put(5,0){\line(0,1){20}}\put(10,0){\line(0,1){25}}
\put(15,0){\line(0,1){25}}\put(20,0){\line(0,1){25}}\put(25,5){\line(0,1){15}}
\put(30,10){\line(0,1){10}}\put(35,10){\line(0,1){10}}
%%%%%%%%%%%%%%%%%%%%%%%%%%%%%%%%%%%%%%%%REMPLISSAGE
\put(15,20){\makebox(6,4)[c]{\small \1}}
\put(10,-7){\makebox(6,4)[c]{\small $\c^{(1)}=(1,0)$}}
%%%%%%%%%%%%%%%%%%%%%%%%%%%%%%%%%%%%%%%%
\linethickness{0.5mm} \put(10,20){\line(1,0){10}}
\linethickness{0.7mm}
\put(10,0){\line(1,0){10}}\put(10,25){\line(1,0){10}}
\put(10,0){\line(0,1){25}}\put(20,0){\line(0,1){25}}
\end{picture}}
\hspace{1cm}
%%%%%%%%%%%%%%%%%%%%%%%%%%%%%%%%%%%%%%%%%%%%%%%%%%%%%%%%%%%%%%%%%%%%%%%%%%%%%%%%%%%%%%%%%%%%%%%%STEP2
 {\setlength{\unitlength}{1mm}
\begin{picture}(35,25)(0,0)
\put(-2,3){\makebox(6,4)[c]{$F_2$}}
%%%%%%%%%%%%%%%%%%%%%%%%%%%%%%%%%%%%%%%%coloriage
\put(5,0){\color{red}{\rule{5mm}{20mm}}}
\put(15,0){\color{red}{\rule{5mm}{20mm}}}
\put(10,5){\color{red}{\rule{5mm}{15mm}}}
%%%%%%%%%%%%%%%%%%%%%%%%%%%%%%%%%%%%%%%%DIAGRAMME
\put(5,0){\line(1,0){15}}\put(5,5){\line(1,0){20}}\put(0,10){\line(1,0){35}}
\put(0,15){\line(1,0){35}}\put(0,20){\line(1,0){35}}\put(10,25){\line(1,0){10}}
\put(0,10){\line(0,1){10}}\put(5,0){\line(0,1){20}}\put(10,0){\line(0,1){25}}
\put(15,0){\line(0,1){25}}\put(20,0){\line(0,1){25}}\put(25,5){\line(0,1){15}}
\put(30,10){\line(0,1){10}}\put(35,10){\line(0,1){10}}
%%%%%%%%%%%%%%%%%%%%%%%%%%%%%%%%%%%%%%%%REMPLISSAGE
\put(15,20){\makebox(6,4)[c]{\small \1}}
\put(10,0){\makebox(6,4)[c]{\small \1}}
\put(10,-7){\makebox(6,4)[c]{\small $\c^{(5)}=(0,0)$}}
%%%%%%%%%%%%%%%%%%%%%%%%%%%%%%%%%%%%%%%%
\linethickness{0.5mm} \put(5,5){\line(1,0){15}}
\linethickness{0.7mm}
\put(5,0){\line(1,0){15}}\put(5,20){\line(1,0){15}}
\put(5,0){\line(0,1){20}}\put(20,0){\line(0,1){20}}
\end{picture}}
\end{center}
\vspace{1cm}
\begin{center}
%%%%%%%%%%%%%%%%%%%%%%%%%%%%%%%%%%%%%%%%%%%%%%%%%%%%%%%%%%%%%%%%%%%%%%%%%%%%%%%%%%%%%%%%%%%%%%%%STEP3
 {\setlength{\unitlength}{1mm}
\begin{picture}(35,25)(0,0)
\put(-2,3){\makebox(6,4)[c]{$F_3$}}
%%%%%%%%%%%%%%%%%%%%%%%%%%%%%%%%%%%%%%%%coloriage
\put(5,0){\color{red}{\rule{5mm}{20mm}}}
\put(15,0){\color{red}{\rule{5mm}{20mm}}}
\put(10,5){\color{red}{\rule{5mm}{15mm}}}
%%%%%%%%%%%%%%%%%%%%%%%%%%%%%%%%%%%%%%%%DIAGRAMME
\put(5,0){\line(1,0){15}}\put(5,5){\line(1,0){20}}\put(0,10){\line(1,0){35}}
\put(0,15){\line(1,0){35}}\put(0,20){\line(1,0){35}}\put(10,25){\line(1,0){10}}
\put(0,10){\line(0,1){10}}\put(5,0){\line(0,1){20}}\put(10,0){\line(0,1){25}}
\put(15,0){\line(0,1){25}}\put(20,0){\line(0,1){25}}\put(25,5){\line(0,1){15}}
\put(30,10){\line(0,1){10}}\put(35,10){\line(0,1){10}}
%%%%%%%%%%%%%%%%%%%%%%%%%%%%%%%%%%%%%%%%REMPLISSAGE
\put(15,20){\makebox(6,4)[c]{\small \1}}
\put(10,0){\makebox(6,4)[c]{\small \1}}
\put(10,-7){\makebox(6,4)[c]{\small $\c^{(4)}=(0)$}}
%%%%%%%%%%%%%%%%%%%%%%%%%%%%%%%%%%%%%%%%
\linethickness{0.5mm} \put(5,10){\line(1,0){20}}
\linethickness{0.7mm}
\put(5,5){\line(1,0){20}}\put(5,20){\line(1,0){20}}
\put(5,5){\line(0,1){15}}\put(25,5){\line(0,1){15}}
\end{picture}}
\hspace{1cm}
%%%%%%%%%%%%%%%%%%%%%%%%%%%%%%%%%%%%%%%%%%%%%%%%%%%%%%%%%%%%%%%%%%%%%%%%%%%%%%%%%%%%%%%%%%%%%%%%STEP4
 {\setlength{\unitlength}{1mm}
\begin{picture}(35,25)(0,0)
\put(-2,3){\makebox(6,4)[c]{$F_4$}}
%%%%%%%%%%%%%%%%%%%%%%%%%%%%%%%%%%%%%%%%coloriage
\put(5,0){\color{red}{\rule{5mm}{20mm}}}
\put(15,0){\color{red}{\rule{5mm}{20mm}}}
\put(10,5){\color{red}{\rule{5mm}{15mm}}}
\put(20,10){\color{red}{\rule{5mm}{5mm}}}
\put(25,10){\color{red}{\rule{5mm}{5mm}}}
%%%%%%%%%%%%%%%%%%%%%%%%%%%%%%%%%%%%%%%%DIAGRAMME
\put(5,0){\line(1,0){15}}\put(5,5){\line(1,0){20}}\put(0,10){\line(1,0){35}}
\put(0,15){\line(1,0){35}}\put(0,20){\line(1,0){35}}\put(10,25){\line(1,0){10}}
\put(0,10){\line(0,1){10}}\put(5,0){\line(0,1){20}}\put(10,0){\line(0,1){25}}
\put(15,0){\line(0,1){25}}\put(20,0){\line(0,1){25}}\put(25,5){\line(0,1){15}}
\put(30,10){\line(0,1){10}}\put(35,10){\line(0,1){10}}
%%%%%%%%%%%%%%%%%%%%%%%%%%%%%%%%%%%%%%%%REMPLISSAGE
\put(15,20){\makebox(6,4)[c]{\small \1}}
\put(10,0){\makebox(6,4)[c]{\small \1}}
\put(20,15){\makebox(6,4)[c]{\small \1}}
\put(25,15){\makebox(6,4)[c]{\small \1}}
\put(10,-7){\makebox(6,4)[c]{\small $\c^{(2)}=(1,0,1)$}}
%%%%%%%%%%%%%%%%%%%%%%%%%%%%%%%%%%%%%%%%
\linethickness{0.5mm} \put(0,15){\line(1,0){35}}
\linethickness{0.7mm}
\put(0,10){\line(1,0){35}}\put(0,20){\line(1,0){35}}
\put(0,10){\line(0,1){10}}\put(35,10){\line(0,1){10}}
\end{picture}}
\hspace{1cm}
%%%%%%%%%%%%%%%%%%%%%%%%%%%%%%%%%%%%%%%%%%%%%%%%%%%%%%%%%%%%%%%%%%%%%%%%%%%%%%%%%%%%%%%%%%%%%%%%STEP5
 {\setlength{\unitlength}{1mm}
\begin{picture}(35,25)(0,0)
\put(-2,3){\makebox(6,4)[c]{$F_5$}}
%%%%%%%%%%%%%%%%%%%%%%%%%%%%%%%%%%%%%%%%coloriage
\put(5,0){\color{red}{\rule{5mm}{20mm}}}
\put(15,0){\color{red}{\rule{5mm}{20mm}}}
\put(10,5){\color{red}{\rule{5mm}{15mm}}}
\put(20,10){\color{red}{\rule{5mm}{5mm}}}
\put(25,10){\color{red}{\rule{5mm}{5mm}}}
%%%%%%%%%%%%%%%%%%%%%%%%%%%%%%%%%%%%%%%%DIAGRAMME
\put(5,0){\line(1,0){15}}\put(5,5){\line(1,0){20}}\put(0,10){\line(1,0){35}}
\put(0,15){\line(1,0){35}}\put(0,20){\line(1,0){35}}\put(10,25){\line(1,0){10}}
\put(0,10){\line(0,1){10}}\put(5,0){\line(0,1){20}}\put(10,0){\line(0,1){25}}
\put(15,0){\line(0,1){25}}\put(20,0){\line(0,1){25}}\put(25,5){\line(0,1){15}}
\put(30,10){\line(0,1){10}}\put(35,10){\line(0,1){10}}
%%%%%%%%%%%%%%%%%%%%%%%%%%%%%%%%%%%%%%%%REMPLISSAGE
\put(15,20){\makebox(6,4)[c]{\small \1}}
\put(10,0){\makebox(6,4)[c]{\small \1}}
\put(20,15){\makebox(6,4)[c]{\small \1}}
\put(25,15){\makebox(6,4)[c]{\small \1}}
\put(0,10){\makebox(6,4)[c]{\small \1}}
\put(30,10){\makebox(6,4)[c]{\small \1}}
\put(10,-7){\makebox(6,4)[c]{\small $\c^{(3)}=(0,0,0)$}}
%%%%%%%%%%%%%%%%%%%%%%%%%%%%%%%%%%%%%%%%
\linethickness{0.7mm}
\put(0,10){\line(1,0){35}}\put(0,15){\line(1,0){35}}
\put(0,10){\line(0,1){5}}\put(35,10){\line(0,1){5}}
\end{picture}}
\vspace{0.5cm}
 \caption{ The step-by-step construction of $\Upsilon({\rm{\textbf{c}}})$}\label{fig:Upsilon}
\end{center}
\end{figure}
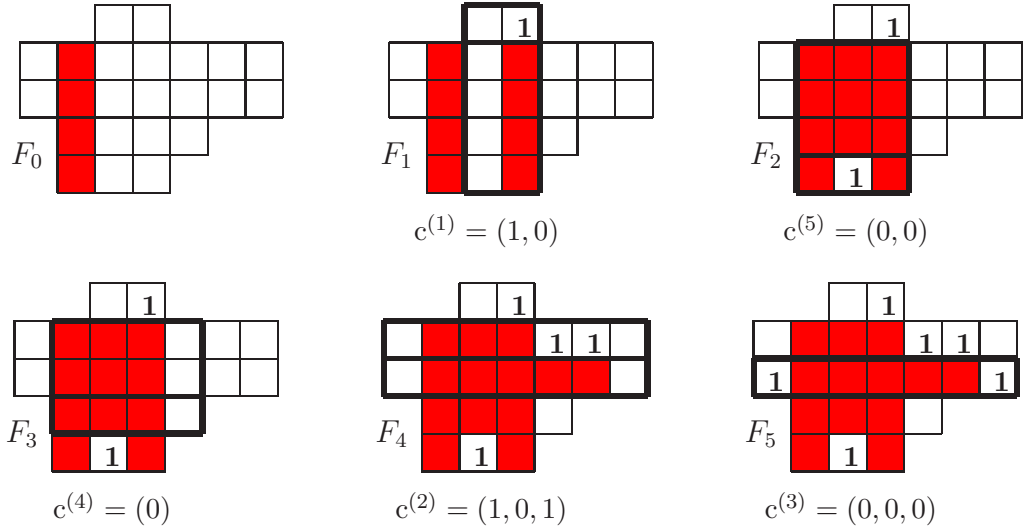

Suppose
${\textrm{\textbf{c}}}=(\c^{(1)},\c^{(2)},\c^{(3)},\c^{(4)},\c^{(5)})$
with $\c^{(1)}=(1,0)$, $\c^{(2)}=(1,0,1)$, $\c^{(3)}=(0,0,0)$,
$\c^{(4)}=(0)$ and $\c^{(5)}=(0,0)$.  The step by step construction
of $\Upsilon({\textrm{\textbf{c}}})$ is given in
Figure~\ref{fig:Upsilon}.
It is easily seen that $\Upsilon$ is the reverse of $\Psi$, and thus
$\Psi$ is bijective. Now let
${\textrm{\textbf{c}}}=(\c^{(1)},\c^{(2)},\ldots,\c^{(s)})\in\C_{m_{1}+1}(h_{1}-m_1)\times\C_{m_{2}+1}(h_{2}-m_2)\times
\cdots\times \C_{m_{s}+1}(h_{s}-m_s)$ and $F=\Upsilon({\bf c})$ be
the corresponding $01$-filling. Let $i$ be an integer in $[s]$ and
$ce$ be the cell of the $i$-th row $R_i$ of $F$ which contains the
$j$-th $1$ of $R_i$ (from  left to right as usual). It then follows
from the definition of $\Upsilon$ that
\begin{align*}
{\luc}(ce;F)=\c_{1}^{(i)}+\c_{2}^{(i)}+\cdots+\c_{j}^{(i)}\quad\text{and}\quad
{\ruc}(ce;F)=\c_{j+1}^{(i)}+\c_{j+2}^{(i)}+\cdots+\c_{m_i+1}^{(i)}.
\end{align*}

We summarize the main properties of $\Psi$ in Theorem~\ref{thm:Psi}.
\begin{thm}\label{thm:Psi}
The map $\Psi:\N^c(T,{\bf m};A)\to \C_{m_{1}+1}(h_{1}-m_1)\times
\cdots\times \C_{m_{s}+1}(h_{s}-m_s)$ is a bijection such that for
any $F\in\N^c(T,\textbf{m};A)$ and any cell $ce$ in $F$, if
$\Psi(F)=(\c^{(1)},\c^{(2)},\cdots,\c^{(s)})$, we have
\begin{align*}
{\luc}(ce;F)&=\c_{1}^{(i)}+\c_{2}^{(i)}+\cdots+\c_{j}^{(i)}\\
{\ruc}(ce;F)&=\c_{j+1}^{(i)}+\c_{j+2}^{(i)}+\cdots+\c_{m_i+1}^{(i)}.
\end{align*}
\end{thm}

\subsection{Proof of~\eqref{eq:symfilling}}
  The proof is based on the correspondence $\Psi$ and the following
identity, which can be easily proved (by induction for instance).

\begin{lem}\label{lem:pq-id}
For any integers $n\geq k\geq0$,
$$
\sum_{(\c_1,\c_2,\ldots,\c_{k+1})\in\C_{k+1}(n)}\prod_{j=1}^k
p^{\sum_{r=1}^j \c_{r}}\,q^{\sum_{r=j+1}^{k+1}\c_{r}}={n+k\brack
k}_{p,q}.
$$
\end{lem}

Suppose $Up(T)=\{R_1,R_2,\ldots,R_{i_0}\}$ and
$Low(T)=\{R_{i_0+1},\ldots,R_s\}$. By \eqref{eq:cro-sum} and
\eqref{eq:ne-sum} we have

\begin{align*}
\sum_{F\in\N^c(T,\textbf{m};A)}p^{{\nec}_2(F)}q^{{\sec}_2(F)}&=
\sum_{F\in\N^c(T,\textbf{m};A)}\prod_{ce\in
Up(T)}p^{{\luc}(ce;F)}q^{{\ruc}(ce;F)}\prod_{ce\in
Low(T)}p^{{\ruc}(ce;F)}q^{{\luc}(ce;F)}\\
&= \sum_{F\in\N^c(T,\textbf{m};A)}\prod_{i=1}^{i_0}\prod_{ce\in
R_i}p^{{\luc}(ce;F)}q^{{\ruc}(ce;F)}\prod_{i=i_0+1}^{s}\prod_{ce\in
R_i}p^{{\ruc}(ce;F)}q^{{\luc}(ce;F)}.
\end{align*}

Let ${\bf \C}=\C_{m_{1}+1}(h_{1}-m_1)\times \cdots\times
\C_{m_{s}+1}(h_{s}-m_s)$. It follows from Theorem~\ref{thm:Psi} that
the right-hand side of the last equality can be rewritten
\begin{align*}
&\sum_{(\c^{(1)},\c^{(2)},\ldots,\c^{(s)})\in{\bf \C}}
 \prod_{i=1}^{i_0}\left(\prod_{j=1}^{m_i}p^{\sum_{r=1}^j
\c^{(i)}_{r}}\,q^{\sum_{r=j+1}^{m_i+1}\c^{(i)}_{r}}\right)
 \prod_{i=i_0+1}^{s}\left(\prod_{j=1}^{m_i} p^{\sum_{r=j+1}^{m_i+1}\c^{(i)}_{r}}\,q^{\sum_{r=1}^j
\c^{(i)}_{r}}\right)\\
&=\prod_{i=1}^{i_0}\left( \sum_{\c^{(i)}\in\C_{m_i+1}(h_i-m_i)}
 \prod_{j=1}^{m_i}p^{\sum_{r=1}^j
\c^{(i)}_{r}}\,q^{\sum_{r=j+1}^{m_i+1}\c^{(i)}_{r}}\right)
 \prod_{i=i_0+1}^{s}
\left( \sum_{\c^{(i)}\in\C_{m_i+1}(h_i-m_i)}\prod_{j=1}^{m_i}
p^{\sum_{r=j+1}^{m_i+1}\c^{(i)}_{r}}\,q^{\sum_{r=1}^j
\c^{(i)}_{r}}\right).\\
\end{align*}

 Applying Lemma~\ref{lem:pq-id} conclude the proof of~\eqref{eq:symfilling}
and thus of Theorem~\ref{thm:distribution}.

%%%%****%%%%****%%%%****%%%%****%%%%****%%%%****%%%%****%%%%*
%%%%
%%%% New Section
%%%%
%%%%****%%%%****%%%%****%%%%****%%%%****%%%%****%%%%****%%%%*

\section{A bijective proof of Corollary~\ref{thm:refinement main1}}

 Let $T$ be a moon polyomino with $s$~rows and $t$ columns.
In this section, we present a mapping $\Phi$ such that for any
$s$-uple of positive integers ${\bf m}=(m_1,m_2,\ldots,m_s)$ and set
$A$ of positive integers, the map $\Phi$ is a bijection
$\Phi:\N^c(T,{\bf m};A)\to \N^c(T,{\bf m};A)$ such that for any
$F\in \N^c(T,{\bf m};A)$, we have
\begin{align*}
({\nec}_2,{\sec}_2)(\Phi(F))=({\sec}_2,{\nec}_2)(F).
\end{align*}
This gives a direct combinatorial proof of the symmetry of the joint
distribution of $({\nec}_2,{\sec}_2)$ over each $\N^c(T,{\bf m};A)$,
$\N(T;A,B)$ (set $m_i=0$ for $i\in B$ and $1$ otherwise), and
$\N^r(T,{\bf n};B)$ (compose with the rotation about
$90^{\circ}$).\\

In fact, the map $\Phi$ is just a byproduct of the constructions
given in the previous section. It is also a generalization of an
involution presented in~\cite{KaZe} to prove the symmetry of
$({\cro}_2,{\ne}_2)$ over set partitions and matchings.

Let $\c=(\c_1,\c_2,\ldots,\c_k)$ be a composition. Define the
\emph{reverse} ${\rev}(\c)$ of $\c$ as the composition
${\rev}(\c)=(\c_k,\c_{k-1},\ldots,\c_1)$. Given a sequence of
compositions ${\bf \c}=(\c^{(1)},\c^{(2)},\ldots,\c^{(s)})$, we set
$\Rev({\bf
\c})=(\rev(\c^{(1)}),\rev(\c^{(2)}),\ldots,\rev(\c^{(s)}))$.

Let $\Phi:\N^c(T,\textbf{m};A)\mapsto \N^c(T,\textbf{m};A)$ be the
map defined by
$$
\Phi=\Psi^{-1}\circ\Rev\circ \Psi=\Upsilon\circ\Rev\circ \Psi.
$$

The following proposition is an immediate consequence of the
properties of $\Psi$ (see Theorem~\ref{thm:Psi}).

\begin{prop}
The map $\Phi$ is an involution on $\N^c(T,\mathbf{m};A)$ such that
for any $F\in \N^c(T,\mathbf{m};A)$ and any cell $c$ of $T$, we have
\begin{align*}
{\luc}(c;\Phi(F))={\ruc}(c;F)\quad\text{and}\quad
{\ruc}(c;\Phi(F))={\luc}(c;F).
\end{align*}
 In particular, we have $({\nec}_2,{\sec}_2)(\Phi(F))=({\sec}_2,{\nec}_2)(F)$.
\end{prop}

It could be useful to give a direct description of $\Phi$. Let $F\in
\N^c(T,{\bf m};A)$.

 (1) Color the columns of polyomino $T$ indexed by the set $A$. Denote by $F'_0$ the
 colored polyomino obtained.

(2) Contruct a sequence of colored fillings $(F'_j)_{j=1...s}$ of
$T$ as follows. Suppose $R_{i_1}\prec R_{i_2} \prec \cdots \prec
R_{i_s}$. Then for $j$ from 1 to $s$, the (colored) filling $F'_{j}$
is obtained from $F'_{j-1}$ as follows:
\begin{itemize}
\item if $m_{i_j}=0$, do nothing,
\item else, read the $m_{i_j}$ 1's in the ${i_j}$-th row of
$F$ from left to right and denote the number of uncolored cells (in
the coloring of $F$) strictly
\begin{itemize}
\item to the left of the first~$1$ by $t_{0}$,
\item between the $u$-th 1 and the $(u+1)$-th 1, $1\leq u\leq m_{i_j}-1$, by $t_{u}$,
\item to the right of the last~$1$ by $t_{m_{i_j}}$.
\end{itemize}
 Then insert $m_{i_j}$ 1's in the $i_j\,$-th row of $F'_{j-1}$ in such a
way that the number of uncolored cells on this row strictly
\begin{itemize}
\item to the left of the first~$1$ is $t_{m_{i_j}}$,
\item between the $u$-th 1 and the $(u+1)$-th 1, $1\leq u\leq m_{i_j}-1$, is $t_{m_{i_j}-u}$,
\item to the right of the last~$1$ is $t_{0}$.
\end{itemize}
Next, color the cells which are both contained in the $i_j$-th
rectangle and below (resp., above) the new $1$'s inserted in
$R_{i_j}$ if $R_{i_j}\in Up(T)$ (resp., $R_{i_j}\in Low(T)$).
\end{itemize}

(3) Set $\Phi(F)=F'_s$. For a better understanding, we give an
illustration. Suppose $F$ is the filling given below. \begin{center}
{\setlength{\unitlength}{0.8mm}
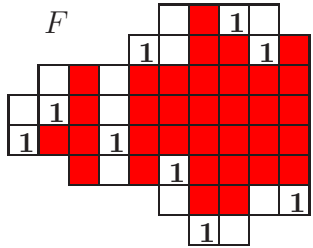
\begin{picture}(50,40)(0,0)
\put(5,35){\makebox(6,4)[c]{$F$}}
%%%%%%%%%%%%%%%%%%%%%%%%%%%%%%%%%%%%%%%%coloriage
\put(30,5){\color{red}{\rule{4mm}{28mm}}}\put(35,5){\color{red}{\rule{4mm}{24mm}}}
\put(45,10){\color{red}{\rule{4mm}{20mm}}}\put(25,15){\color{red}{\rule{4mm}{12mm}}}
\put(5,15){\color{red}{\rule{4mm}{4mm}}}\put(10,10){\color{red}{\rule{4mm}{16mm}}}
\put(20,10){\color{red}{\rule{4mm}{16mm}}}\put(40,10){\color{red}{\rule{4mm}{16mm}}}
%%%%%%%%%%%%%%%%%%%%%%%%%%%%%%%%%%%%%%%%DIAGRAMME
\put(30,0){\line(1,0){10}} \put(25,5){\line(1,0){25}}
\put(10,10){\line(1,0){40}} \put(0,15){\line(1,0){50}}
\put(0,20){\line(1,0){50}} \put(0,25){\line(1,0){50}}
\put(5,30){\line(1,0){45}} \put(20,35){\line(1,0){30}}
\put(15,30){\line(1,0){10}} \put(25,40){\line(1,0){20}}
\put(0,15){\line(0,1){10}}\put(5,15){\line(0,1){15}}
\put(10,10){\line(0,1){20}}\put(15,10){\line(0,1){20}}\put(20,10){\line(0,1){25}}
\put(25,5){\line(0,1){35}}\put(30,0){\line(0,1){40}}\put(35,0){\line(0,1){40}}\put(40,0){\line(0,1){40}}
\put(45,5){\line(0,1){35}} \put(50,5){\line(0,1){30}}
%%%%%%%%%%%%%%%%%%%%%%%%%%%%%%%%%%%%%%%%REMPLISSAGE
\put(30,0){\makebox(6,4)[c]{\small \1}}
\put(45,5){\makebox(6,4)[c]{\small \1}}
\put(25,10){\makebox(6,4)[c]{\small \1}}
\put(15,15){\makebox(6,4)[c]{\small\1}}
\put(0,15){\makebox(6,4)[c]{\small \1}}
\put(5,20){\makebox(6,4)[c]{\small\1}}
\put(20,30){\makebox(6,4)[c]{\small\1}}
\put(40,30){\makebox(6,4)[c]{\small \1}}
\put(35,35){\makebox(6,4)[c]{\small \1}}
\end{picture}}\end{center}

Then the step-by-step construction of $\Phi(F)$ goes as follows.

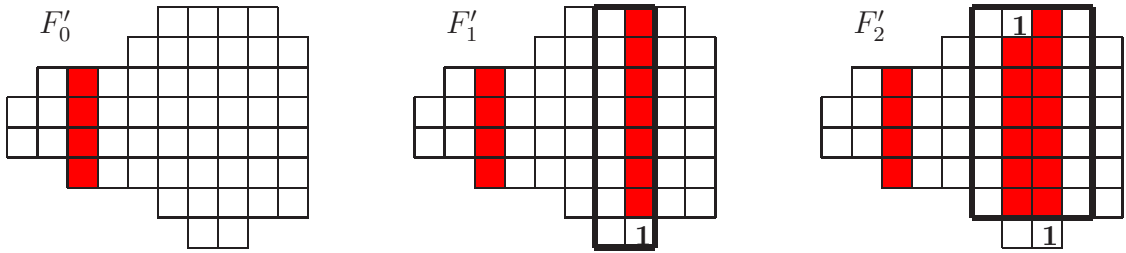
\begin{center}
%%%%%%%%%%%%%%%%%%%%%%%%%%%%%%%%%%%%%%%%%%%%%%%%%%%%%%%%%%%%%%%%%%%%%%%%%%%%%%%%%%%%%%%%%%%%%%%%STEP0
 {\setlength{\unitlength}{0.8mm}
\begin{picture}(50,45)(0,0)
\put(5,35){\makebox(6,4)[c]{$F'_0$}}
%%%%%%%%%%%%%%%%%%%%%%%%%%%%%%%%%%%%%%%%coloriage
\put(10,10){\color{red}{\rule{4mm}{16mm}}}
%%%%%%%%%%%%%%%%%%%%%%%%%%%%%%%%%%%%%%%%DIAGRAMME
\put(30,0){\line(1,0){10}} \put(25,5){\line(1,0){25}}
\put(10,10){\line(1,0){40}} \put(0,15){\line(1,0){50}}
\put(0,20){\line(1,0){50}} \put(0,25){\line(1,0){50}}
\put(5,30){\line(1,0){45}} \put(20,35){\line(1,0){30}}
\put(15,30){\line(1,0){10}} \put(25,40){\line(1,0){20}}
\put(0,15){\line(0,1){10}}\put(5,15){\line(0,1){15}}
\put(10,10){\line(0,1){20}}\put(15,10){\line(0,1){20}}\put(20,10){\line(0,1){25}}
\put(25,5){\line(0,1){35}}\put(30,0){\line(0,1){40}}\put(35,0){\line(0,1){40}}\put(40,0){\line(0,1){40}}
\put(45,5){\line(0,1){35}} \put(50,5){\line(0,1){30}}
\end{picture}}
\hspace{1cm}
%%%%%%%%%%%%%%%%%%%%%%%%%%%%%%%%%%%%%%%%%%%%%%%%%%%%%%%%%%%%%%%%%%%%%%%%%%%%%%%%%%%%%%%%%%%%%%%%STEP1
 {\setlength{\unitlength}{0.8mm}
\begin{picture}(50,45)(0,0)
\put(5,35){\makebox(6,4)[c]{$F'_1$}}
%%%%%%%%%%%%%%%%%%%%%%%%%%%%%%%%%%%%%%%%coloriage
\put(10,10){\color{red}{\rule{4mm}{16mm}}}
\put(35,5){\color{red}{\rule{4mm}{28mm}}}
%%%%%%%%%%%%%%%%%%%%%%%%%%%%%%%%%%%%%%%%DIAGRAMME
\put(30,0){\line(1,0){10}} \put(25,5){\line(1,0){25}}
\put(10,10){\line(1,0){40}} \put(0,15){\line(1,0){50}}
\put(0,20){\line(1,0){50}} \put(0,25){\line(1,0){50}}
\put(5,30){\line(1,0){45}} \put(20,35){\line(1,0){30}}
\put(15,30){\line(1,0){10}} \put(25,40){\line(1,0){20}}
\put(0,15){\line(0,1){10}}\put(5,15){\line(0,1){15}}
\put(10,10){\line(0,1){20}}\put(15,10){\line(0,1){20}}\put(20,10){\line(0,1){25}}
\put(25,5){\line(0,1){35}}\put(30,0){\line(0,1){40}}\put(35,0){\line(0,1){40}}\put(40,0){\line(0,1){40}}
\put(45,5){\line(0,1){35}} \put(50,5){\line(0,1){30}}
%%%%%%%%%%%%%%%%%%%%%%%%%%%%%%%%%%%%%%%%REMPLISSAGE
\put(35,0){\makebox(6,4)[c]{\small \1}}
%%%%%%%%%%%%%%%%%%%%%%%%%%%%%%%%%%%%%%%%
\linethickness{0.7mm}
\put(30,0){\line(1,0){10}}\put(30,40){\line(1,0){10}}
\put(30,0){\line(0,1){40}}\put(40,0){\line(0,1){40}}
\end{picture}}
\hspace{1cm}
%%%%%%%%%%%%%%%%%%%%%%%%%%%%%%%%%%%%%%%%%%%%%%%%%%%%%%%%%%%%%%%%%%%%%%%%%%%%%%%%%%%%%%%%%%%%%%%%STEP2
 {\setlength{\unitlength}{0.8mm}
\begin{picture}(50,45)(0,0)
\put(5,35){\makebox(6,4)[c]{$F'_2$}}
%%%%%%%%%%%%%%%%%%%%%%%%%%%%%%%%%%%%%%%%coloriage
\put(10,10){\color{red}{\rule{4mm}{16mm}}}
\put(35,5){\color{red}{\rule{4mm}{28mm}}}
\put(30,5){\color{red}{\rule{4mm}{24mm}}}
%%%%%%%%%%%%%%%%%%%%%%%%%%%%%%%%%%%%%%%%DIAGRAMME
\put(30,0){\line(1,0){10}} \put(25,5){\line(1,0){25}}
\put(10,10){\line(1,0){40}} \put(0,15){\line(1,0){50}}
\put(0,20){\line(1,0){50}} \put(0,25){\line(1,0){50}}
\put(5,30){\line(1,0){45}} \put(20,35){\line(1,0){30}}
\put(15,30){\line(1,0){10}} \put(25,40){\line(1,0){20}}
\put(0,15){\line(0,1){10}}\put(5,15){\line(0,1){15}}
\put(10,10){\line(0,1){20}}\put(15,10){\line(0,1){20}}\put(20,10){\line(0,1){25}}
\put(25,5){\line(0,1){35}}\put(30,0){\line(0,1){40}}\put(35,0){\line(0,1){40}}\put(40,0){\line(0,1){40}}
\put(45,5){\line(0,1){35}} \put(50,5){\line(0,1){30}}
%%%%%%%%%%%%%%%%%%%%%%%%%%%%%%%%%%%%%%%%REMPLISSAGE
\put(35,0){\makebox(6,4)[c]{\small \1}}
\put(30,35){\makebox(6,4)[c]{\small \1}}
%%%%%%%%%%%%%%%%%%%%%%%%%%%%%%%%%%%%%%%%
\linethickness{0.7mm}
\put(25,5){\line(0,1){35}}\put(45,5){\line(0,1){35}}
\put(25,5){\line(1,0){20}}\put(25,40){\line(1,0){20}}
\end{picture}}
\end{center}

\begin{center}
%%%%%%%%%%%%%%%%%%%%%%%%%%%%%%%%%%%%%%%%%%%%%%%%%%%%%%%%%%%%%%%%%%%%%%%%%%%%%%%%%%%%%%%%%%%%%%%%STEP3
 {\setlength{\unitlength}{0.8mm}
\begin{picture}(50,50)(0,0)
\put(5,35){\makebox(6,4)[c]{$F'_3$}}
%%%%%%%%%%%%%%%%%%%%%%%%%%%%%%%%%%%%%%%%coloriage
\put(10,10){\color{red}{\rule{4mm}{16mm}}}
\put(35,5){\color{red}{\rule{4mm}{28mm}}}
\put(30,5){\color{red}{\rule{4mm}{24mm}}}
\put(25,10){\color{red}{\rule{4mm}{20mm}}}
%%%%%%%%%%%%%%%%%%%%%%%%%%%%%%%%%%%%%%%%DIAGRAMME
\put(30,0){\line(1,0){10}} \put(25,5){\line(1,0){25}}
\put(10,10){\line(1,0){40}} \put(0,15){\line(1,0){50}}
\put(0,20){\line(1,0){50}} \put(0,25){\line(1,0){50}}
\put(5,30){\line(1,0){45}} \put(20,35){\line(1,0){30}}
\put(15,30){\line(1,0){10}} \put(25,40){\line(1,0){20}}
\put(0,15){\line(0,1){10}}\put(5,15){\line(0,1){15}}
\put(10,10){\line(0,1){20}}\put(15,10){\line(0,1){20}}\put(20,10){\line(0,1){25}}
\put(25,5){\line(0,1){35}}\put(30,0){\line(0,1){40}}\put(35,0){\line(0,1){40}}\put(40,0){\line(0,1){40}}
\put(45,5){\line(0,1){35}} \put(50,5){\line(0,1){30}}
%%%%%%%%%%%%%%%%%%%%%%%%%%%%%%%%%%%%%%%%REMPLISSAGE
\put(35,0){\makebox(6,4)[c]{\small \1}}
\put(30,35){\makebox(6,4)[c]{\small \1}}
\put(25,5){\makebox(6,4)[c]{\small \1}}
%%%%%%%%%%%%%%%%%%%%%%%%%%%%%%%%%%%%%%%%
\linethickness{0.7mm}
\put(25,5){\line(0,1){30}}\put(50,5){\line(0,1){30}}
\put(25,5){\line(1,0){25}}\put(25,35){\line(1,0){25}}
\end{picture}}
\hspace{1cm}
%%%%%%%%%%%%%%%%%%%%%%%%%%%%%%%%%%%%%%%%%%%%%%%%%%%%%%%%%%%%%%%%%%%%%%%%%%%%%%%%%%%%%%%%%%%%%%%%STEP4
 {\setlength{\unitlength}{0.8mm}
\begin{picture}(50,50)(0,0)
\put(5,35){\makebox(6,4)[c]{$F'_4$}}
%%%%%%%%%%%%%%%%%%%%%%%%%%%%%%%%%%%%%%%%coloriage
\put(10,10){\color{red}{\rule{4mm}{16mm}}}
\put(35,5){\color{red}{\rule{4mm}{28mm}}}
\put(30,5){\color{red}{\rule{4mm}{24mm}}}
\put(25,10){\color{red}{\rule{4mm}{20mm}}}
\put(20,10){\color{red}{\rule{4mm}{16mm}}}
\put(45,10){\color{red}{\rule{4mm}{16mm}}}
%%%%%%%%%%%%%%%%%%%%%%%%%%%%%%%%%%%%%%%%DIAGRAMME
\put(30,0){\line(1,0){10}} \put(25,5){\line(1,0){25}}
\put(10,10){\line(1,0){40}} \put(0,15){\line(1,0){50}}
\put(0,20){\line(1,0){50}} \put(0,25){\line(1,0){50}}
\put(5,30){\line(1,0){45}} \put(20,35){\line(1,0){30}}
\put(15,30){\line(1,0){10}} \put(25,40){\line(1,0){20}}
\put(0,15){\line(0,1){10}}\put(5,15){\line(0,1){15}}
\put(10,10){\line(0,1){20}}\put(15,10){\line(0,1){20}}\put(20,10){\line(0,1){25}}
\put(25,5){\line(0,1){35}}\put(30,0){\line(0,1){40}}\put(35,0){\line(0,1){40}}\put(40,0){\line(0,1){40}}
\put(45,5){\line(0,1){35}} \put(50,5){\line(0,1){30}}
%%%%%%%%%%%%%%%%%%%%%%%%%%%%%%%%%%%%%%%%REMPLISSAGE
\put(35,0){\makebox(6,4)[c]{\small \1}}
\put(30,35){\makebox(6,4)[c]{\small \1}}
\put(25,5){\makebox(6,4)[c]{\small \1}}
\put(20,30){\makebox(6,4)[c]{\small \1}}
\put(45,30){\makebox(6,4)[c]{\small \1}}
%%%%%%%%%%%%%%%%%%%%%%%%%%%%%%%%%%%%%%%%
\linethickness{0.7mm}
\put(20,10){\line(0,1){25}}\put(50,10){\line(0,1){25}}
\put(20,10){\line(1,0){30}}\put(20,35){\line(1,0){30}}
\end{picture}}
\hspace{1cm}
%%%%%%%%%%%%%%%%%%%%%%%%%%%%%%%%%%%%%%%%%%%%%%%%%%%%%%%%%%%%%%%%%%%%%%%%%%%%%%%%%%%%%%%%%%%%%%%%STEP5
{\setlength{\unitlength}{0.8mm}
\begin{picture}(50,50)(0,0)
\put(5,35){\makebox(6,4)[c]{$F'_5$}}
%%%%%%%%%%%%%%%%%%%%%%%%%%%%%%%%%%%%%%%%coloriage
\put(10,10){\color{red}{\rule{4mm}{16mm}}}
\put(35,5){\color{red}{\rule{4mm}{28mm}}}
\put(30,5){\color{red}{\rule{4mm}{24mm}}}
\put(25,10){\color{red}{\rule{4mm}{20mm}}}
\put(20,10){\color{red}{\rule{4mm}{16mm}}}
\put(45,10){\color{red}{\rule{4mm}{16mm}}}
\put(15,15){\color{red}{\rule{4mm}{12mm}}}
%%%%%%%%%%%%%%%%%%%%%%%%%%%%%%%%%%%%%%%%DIAGRAMME
\put(30,0){\line(1,0){10}} \put(25,5){\line(1,0){25}}
\put(10,10){\line(1,0){40}} \put(0,15){\line(1,0){50}}
\put(0,20){\line(1,0){50}} \put(0,25){\line(1,0){50}}
\put(5,30){\line(1,0){45}} \put(20,35){\line(1,0){30}}
\put(15,30){\line(1,0){10}} \put(25,40){\line(1,0){20}}
\put(0,15){\line(0,1){10}}\put(5,15){\line(0,1){15}}
\put(10,10){\line(0,1){20}}\put(15,10){\line(0,1){20}}\put(20,10){\line(0,1){25}}
\put(25,5){\line(0,1){35}}\put(30,0){\line(0,1){40}}\put(35,0){\line(0,1){40}}\put(40,0){\line(0,1){40}}
\put(45,5){\line(0,1){35}} \put(50,5){\line(0,1){30}}
%%%%%%%%%%%%%%%%%%%%%%%%%%%%%%%%%%%%%%%%REMPLISSAGE
\put(35,0){\makebox(6,4)[c]{\small \1}}
\put(30,35){\makebox(6,4)[c]{\small \1}}
\put(25,5){\makebox(6,4)[c]{\small \1}}
\put(20,30){\makebox(6,4)[c]{\small \1}}
\put(45,30){\makebox(6,4)[c]{\small \1}}
\put(15,10){\makebox(6,4)[c]{\small \1}}
%%%%%%%%%%%%%%%%%%%%%%%%%%%%%%%%%%%%%%%%
\linethickness{0.7mm}
\put(10,10){\line(0,1){20}}\put(50,10){\line(0,1){20}}
\put(10,10){\line(1,0){40}}\put(10,30){\line(1,0){40}}
\end{picture}}
\end{center}

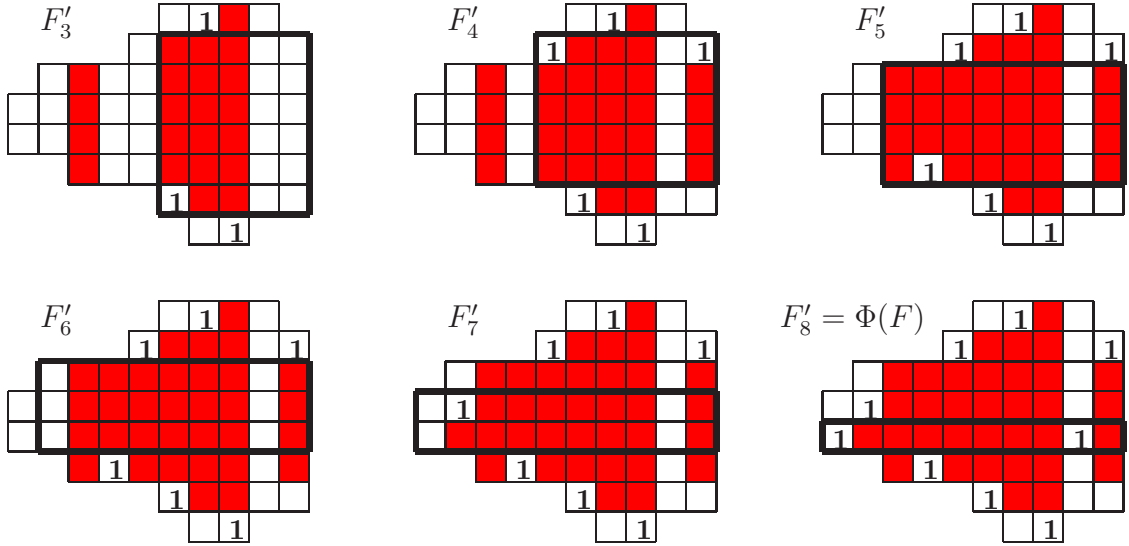
\begin{figure}[h!]\label{fig:invo-fillings}

\vspace{-0.1cm}

\vspace{-0.3cm}

\vspace{-0.3cm}
\begin{center}
%%%%%%%%%%%%%%%%%%%%%%%%%%%%%%%%%%%%%%%%%%%%%%%%%%%%%%%%%%%%%%%%%%%%%%%%%%%%%%%%%%%%%%%%%%%%%%%%STEP6
 {\setlength{\unitlength}{0.8mm}
\begin{picture}(50,50)(0,0)
\put(5,35){\makebox(6,4)[c]{$F'_6$}}
%%%%%%%%%%%%%%%%%%%%%%%%%%%%%%%%%%%%%%%%coloriage
\put(10,10){\color{red}{\rule{4mm}{16mm}}}
\put(35,5){\color{red}{\rule{4mm}{28mm}}}
\put(30,5){\color{red}{\rule{4mm}{24mm}}}
\put(25,10){\color{red}{\rule{4mm}{20mm}}}
\put(20,10){\color{red}{\rule{4mm}{16mm}}}
\put(45,10){\color{red}{\rule{4mm}{16mm}}}
\put(15,15){\color{red}{\rule{4mm}{12mm}}}
%%%%%%%%%%%%%%%%%%%%%%%%%%%%%%%%%%%%%%%%DIAGRAMME
\put(30,0){\line(1,0){10}} \put(25,5){\line(1,0){25}}
\put(10,10){\line(1,0){40}} \put(0,15){\line(1,0){50}}
\put(0,20){\line(1,0){50}} \put(0,25){\line(1,0){50}}
\put(5,30){\line(1,0){45}} \put(20,35){\line(1,0){30}}
\put(15,30){\line(1,0){10}} \put(25,40){\line(1,0){20}}
\put(0,15){\line(0,1){10}}\put(5,15){\line(0,1){15}}
\put(10,10){\line(0,1){20}}\put(15,10){\line(0,1){20}}\put(20,10){\line(0,1){25}}
\put(25,5){\line(0,1){35}}\put(30,0){\line(0,1){40}}\put(35,0){\line(0,1){40}}\put(40,0){\line(0,1){40}}
\put(45,5){\line(0,1){35}} \put(50,5){\line(0,1){30}}
%%%%%%%%%%%%%%%%%%%%%%%%%%%%%%%%%%%%%%%%REMPLISSAGE
\put(35,0){\makebox(6,4)[c]{\small \1}}
\put(30,35){\makebox(6,4)[c]{\small \1}}
\put(25,5){\makebox(6,4)[c]{\small \1}}
\put(20,30){\makebox(6,4)[c]{\small \1}}
\put(45,30){\makebox(6,4)[c]{\small \1}}
\put(15,10){\makebox(6,4)[c]{\small \1}}
%%%%%%%%%%%%%%%%%%%%%%%%%%%%%%%%%%%%%%%%
\linethickness{0.7mm}
\put(5,15){\line(0,1){15}}\put(50,15){\line(0,1){15}}
\put(5,15){\line(1,0){45}}\put(5,30){\line(1,0){45}}
\end{picture}}
\hspace{1cm}
%%%%%%%%%%%%%%%%%%%%%%%%%%%%%%%%%%%%%%%%%%%%%%%%%%%%%%%%%%%%%%%%%%%%%%%%%%%%%%%%%%%%%%%%%%%%%%%%STEP7
{\setlength{\unitlength}{0.8mm}
\begin{picture}(50,50)(0,0)
\put(5,35){\makebox(6,4)[c]{$F'_7$}}
%%%%%%%%%%%%%%%%%%%%%%%%%%%%%%%%%%%%%%%%coloriage
\put(10,10){\color{red}{\rule{4mm}{16mm}}}
\put(35,5){\color{red}{\rule{4mm}{28mm}}}
\put(30,5){\color{red}{\rule{4mm}{24mm}}}
\put(25,10){\color{red}{\rule{4mm}{20mm}}}
\put(20,10){\color{red}{\rule{4mm}{16mm}}}
\put(45,10){\color{red}{\rule{4mm}{16mm}}}
\put(15,15){\color{red}{\rule{4mm}{12mm}}}
\put(5,15){\color{red}{\rule{4mm}{4mm}}}
%%%%%%%%%%%%%%%%%%%%%%%%%%%%%%%%%%%%%%%%DIAGRAMME
\put(30,0){\line(1,0){10}} \put(25,5){\line(1,0){25}}
\put(10,10){\line(1,0){40}} \put(0,15){\line(1,0){50}}
\put(0,20){\line(1,0){50}} \put(0,25){\line(1,0){50}}
\put(5,30){\line(1,0){45}} \put(20,35){\line(1,0){30}}
\put(15,30){\line(1,0){10}} \put(25,40){\line(1,0){20}}
\put(0,15){\line(0,1){10}}\put(5,15){\line(0,1){15}}
\put(10,10){\line(0,1){20}}\put(15,10){\line(0,1){20}}\put(20,10){\line(0,1){25}}
\put(25,5){\line(0,1){35}}\put(30,0){\line(0,1){40}}\put(35,0){\line(0,1){40}}\put(40,0){\line(0,1){40}}
\put(45,5){\line(0,1){35}} \put(50,5){\line(0,1){30}}
%%%%%%%%%%%%%%%%%%%%%%%%%%%%%%%%%%%%%%%%REMPLISSAGE
\put(35,0){\makebox(6,4)[c]{\small \1}}
\put(30,35){\makebox(6,4)[c]{\small \1}}
\put(25,5){\makebox(6,4)[c]{\small \1}}
\put(20,30){\makebox(6,4)[c]{\small \1}}
\put(45,30){\makebox(6,4)[c]{\small \1}}
\put(15,10){\makebox(6,4)[c]{\small \1}}
\put(5,20){\makebox(6,4)[c]{\small \1}}
%%%%%%%%%%%%%%%%%%%%%%%%%%%%%%%%%%%%%%%%
\linethickness{0.7mm}
\put(0,15){\line(0,1){10}}\put(50,15){\line(0,1){10}}
\put(0,15){\line(1,0){50}}\put(0,25){\line(1,0){50}}
\end{picture}}
\hspace{1cm}
%%%%%%%%%%%%%%%%%%%%%%%%%%%%%%%%%%%%%%%%%%%%%%%%%%%%%%%%%%%%%%%%%%%%%%%%%%%%%%%%%%%%%%%%%%%%%%%%STEP8
{\setlength{\unitlength}{0.8mm}
\begin{picture}(50,50)(0,0)
\put(2,35){\makebox(6,4)[c]{$F'_8=\Phi(F)$}}
%%%%%%%%%%%%%%%%%%%%%%%%%%%%%%%%%%%%%%%%coloriage
\put(10,10){\color{red}{\rule{4mm}{16mm}}}
\put(35,5){\color{red}{\rule{4mm}{28mm}}}
\put(30,5){\color{red}{\rule{4mm}{24mm}}}
\put(25,10){\color{red}{\rule{4mm}{20mm}}}
\put(20,10){\color{red}{\rule{4mm}{16mm}}}
\put(45,10){\color{red}{\rule{4mm}{16mm}}}
\put(15,15){\color{red}{\rule{4mm}{12mm}}}
\put(5,15){\color{red}{\rule{4mm}{4mm}}}
%%%%%%%%%%%%%%%%%%%%%%%%%%%%%%%%%%%%%%%%DIAGRAMME
\put(30,0){\line(1,0){10}} \put(25,5){\line(1,0){25}}
\put(10,10){\line(1,0){40}} \put(0,15){\line(1,0){50}}
\put(0,20){\line(1,0){50}} \put(0,25){\line(1,0){50}}
\put(5,30){\line(1,0){45}} \put(20,35){\line(1,0){30}}
\put(15,30){\line(1,0){10}} \put(25,40){\line(1,0){20}}
\put(0,15){\line(0,1){10}}\put(5,15){\line(0,1){15}}
\put(10,10){\line(0,1){20}}\put(15,10){\line(0,1){20}}\put(20,10){\line(0,1){25}}
\put(25,5){\line(0,1){35}}\put(30,0){\line(0,1){40}}\put(35,0){\line(0,1){40}}\put(40,0){\line(0,1){40}}
\put(45,5){\line(0,1){35}} \put(50,5){\line(0,1){30}}
%%%%%%%%%%%%%%%%%%%%%%%%%%%%%%%%%%%%%%%%REMPLISSAGE
\put(35,0){\makebox(6,4)[c]{\small \1}}
\put(30,35){\makebox(6,4)[c]{\small \1}}
\put(25,5){\makebox(6,4)[c]{\small \1}}
\put(20,30){\makebox(6,4)[c]{\small \1}}
\put(45,30){\makebox(6,4)[c]{\small \1}}
\put(15,10){\makebox(6,4)[c]{\small \1}}
\put(5,20){\makebox(6,4)[c]{\small \1}}
\put(0,15){\makebox(6,4)[c]{\small \1}}
\put(40,15){\makebox(6,4)[c]{\small \1}}
%%%%%%%%%%%%%%%%%%%%%%%%%%%%%%%%%%%%%%%%
\linethickness{0.7mm}
\put(0,15){\line(0,1){5}}\put(50,15){\line(0,1){5}}
\put(0,15){\line(1,0){50}}\put(0,20){\line(1,0){50}}
\end{picture}}
\end{center}
\caption{ The step-by-step construction of $\Phi(F)$}
\end{figure}

%%%%****%%%%****%%%%****%%%%****%%%%****%%%%****%%%%****%%%%*
%%%%
%%%% New Section
%%%%
%%%%****%%%%****%%%%****%%%%****%%%%****%%%%****%%%%****%%%%*

\section{Concluding remarks}

 It is natural, in view of the results obtained in this paper, to ask  if
the joint distribution of the statistic $({\nec}_2,{\sec}_2)$ is
symmetric over arbitrary $01$-fillings of moon polyominoes, i.e.,
there are no restrictions on the number of $1$'s in columns and
rows. The answer is no by means of the following result. Given a
moon polyomino $T$, recall that  $\N^{01}(T)$ is the set of all
${01}$-fillings of $T$.

\begin{prop}\label{prop:nonsym-graphe}
For any $n\geq5$ the numbers of arbitrary $01$-fillings
of~$\Delta_n$
\begin{itemize}
\item with exactly ${n\choose 4}$ descents is equal to $2^n$,
\item with exactly ${n\choose 4}$ ascents is equal to $16$.
\end{itemize}
In particular, for any $n\geq 5$, the statistics ${\nec}_2$ and
${\sec}_2$ are not equidistributed over $\N^{01}(\Delta_n)$, and
thus the joint distribution of $({\nec}_2,{\sec}_2)$ over
$\N^{01}(\Delta_n)$ is not symmetric.

This also implies that the statistics ${\cro}_2$ and ${\ne}_2$ are
not equidistributed over all simple graphs of $[n]$.
\end{prop}

\pf We give the proof for $n=5,6$ since the reasoning can be
generalized for arbitrary~$n$. Suppose $n=5$. Then one can check
that the arbitrary $01$-fillings of $\Delta_5$ with exactly $5$
descents and those with exactly $5$ ascents have respectively the
following "form"
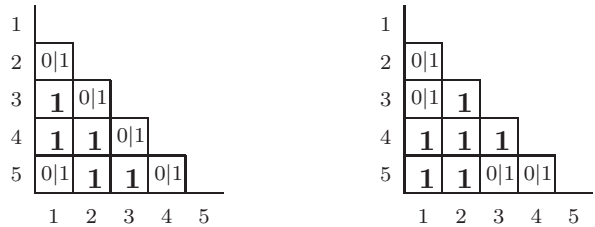
\begin{center}
{\setlength{\unitlength}{1mm}
\begin{picture}(25,30)(0,-5)
%%%%%%%%%%%%%%%%%%%%%%%%%%%%%%%%%%%%%%%%DIAGRAMME
\put(0,0){\line(1,0){25}}\put(0,5){\line(1,0){20}}\put(0,10){\line(1,0){15}}
\put(0,15){\line(1,0){10}}\put(0,20){\line(1,0){5}}
\put(0,0){\line(0,1){25}}\put(5,0){\line(0,1){20}}\put(10,0){\line(0,1){15}}
\put(15,0){\line(0,1){10}}\put(20,0){\line(0,1){5}}
%%%%%%%%%%%%%%%%%%%%%%%%%%%%%%%%%%%%%%%%REMPLISSAGE
\put(0,5){\makebox(6,4)[c]{\small \1}}
\put(0,10){\makebox(6,4)[c]{\small \1}}
\put(5,0){\makebox(6,4)[c]{\small \1}}
\put(5,5){\makebox(6,4)[c]{\small \1}}
\put(10,0){\makebox(6,4)[c]{\small \1}}
\put(0,0){\makebox(5.5,5)[c]{\tiny $0|1$}}
\put(0,15){\makebox(5.5,5)[c]{\tiny $0|1$}}
\put(5,10){\makebox(5.5,5)[c]{\tiny $0|1$}}
\put(10,5){\makebox(5.5,5)[c]{\tiny $0|1$}}
\put(15,0){\makebox(5.5,5)[c]{\tiny $0|1$}}
%%%%%%%%%%%%%%%%%%%%%%%%%%%%%%%%%%%%%%%%ETIQUETAGE DES COLONNES ET LIGNES
\put(0,0){\makebox(5,-6)[c]{\tiny $1$}}
\put(5,0){\makebox(5,-6)[c]{\tiny $2$}}
\put(10,0){\makebox(5,-6)[c]{\tiny $3$}}
\put(15,0){\makebox(5,-6)[c]{\tiny $4$}}
\put(20,0){\makebox(5,-6)[c]{\tiny $5$}}
%%%%%%%%%%%%
\put(0,0){\makebox(-5,5)[c]{\tiny $5$}}
\put(0,5){\makebox(-5,5)[c]{\tiny $4$}}
\put(0,10){\makebox(-5,5)[c]{\tiny $3$}}
\put(0,15){\makebox(-5,5)[c]{\tiny $2$}}
\put(0,20){\makebox(-5,5)[c]{\tiny $1$}}
\end{picture}}
\hspace{2cm} {\setlength{\unitlength}{1mm}
\begin{picture}(25,30)(0,-5)
%%%%%%%%%%%%%%%%%%%%%%%%%%%%%%%%%%%%%%%%DIAGRAMME
\put(0,0){\line(1,0){25}}\put(0,5){\line(1,0){20}}\put(0,10){\line(1,0){15}}
\put(0,15){\line(1,0){10}}\put(0,20){\line(1,0){5}}
\put(0,0){\line(0,1){25}}\put(5,0){\line(0,1){20}}\put(10,0){\line(0,1){15}}
\put(15,0){\line(0,1){10}}\put(20,0){\line(0,1){5}}
%%%%%%%%%%%%%%%%%%%%%%%%%%%%%%%%%%%%%%%%REMPLISSAGE
\put(0,0){\makebox(6,4)[c]{\small \1}}
\put(0,5){\makebox(6,4)[c]{\small \1}}
\put(5,0){\makebox(6,4)[c]{\small \1}}
\put(5,5){\makebox(6,4)[c]{\small \1}}
\put(5,10){\makebox(6,4)[c]{\small \1}}
\put(10,5){\makebox(6,4)[c]{\small \1}}
\put(0,10){\makebox(5.5,5)[c]{\tiny $0|1$}}
\put(0,15){\makebox(5.5,5)[c]{\tiny $0|1$}}
\put(10,0){\makebox(5.5,5)[c]{\tiny $0|1$}}
\put(15,0){\makebox(5.5,5)[c]{\tiny $0|1$}}
%%%%%%%%%%%%%%%%%%%%%%%%%%%%%%%%%%%%%%%%ETIQUETAGE DES COLONNES ET LIGNES
\put(0,0){\makebox(5,-6)[c]{\tiny $1$}}
\put(5,0){\makebox(5,-6)[c]{\tiny $2$}}
\put(10,0){\makebox(5,-6)[c]{\tiny $3$}}
\put(15,0){\makebox(5,-6)[c]{\tiny $4$}}
\put(20,0){\makebox(5,-6)[c]{\tiny $5$}}
%%%%%%%%%%%%
\put(0,0){\makebox(-5,5)[c]{\tiny $5$}}
\put(0,5){\makebox(-5,5)[c]{\tiny $4$}}
\put(0,10){\makebox(-5,5)[c]{\tiny $3$}}
\put(0,15){\makebox(-5,5)[c]{\tiny $2$}}
\put(0,20){\makebox(-5,5)[c]{\tiny $1$}}
\end{picture}}
\end{center}
from which it is easy to obtain the result. Similarly, for $n=6$,
the arbitrary $01$-fillings of $\Delta_6$ with exactly $15$ descents
and those with exactly $15$ ascents have respectively the following
"form"

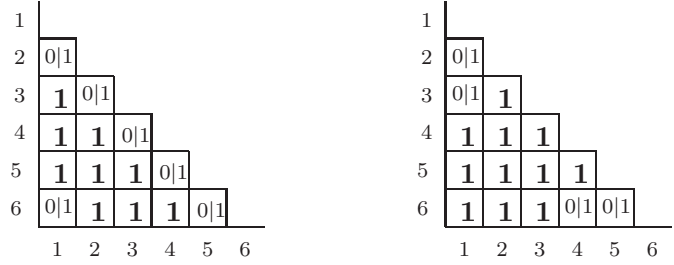
\begin{center}
{\setlength{\unitlength}{1mm}
\begin{picture}(30,35)(0,-5)
%%%%%%%%%%%%%%%%%%%%%%%%%%%%%%%%%%%%%%%%DIAGRAMME
\put(0,0){\line(1,0){30}}\put(0,5){\line(1,0){25}}\put(0,10){\line(1,0){20}}
\put(0,15){\line(1,0){15}}\put(0,20){\line(1,0){10}}\put(0,25){\line(1,0){5}}
\put(0,0){\line(0,1){30}}\put(5,0){\line(0,1){25}}\put(10,0){\line(0,1){20}}
\put(15,0){\line(0,1){15}}\put(20,0){\line(0,1){10}}\put(25,0){\line(0,1){5}}
%%%%%%%%%%%%%%%%%%%%%%%%%%%%%%%%%%%%%%%%REMPLISSAGE
\put(0,0){\makebox(5.5,5)[c]{\tiny $0|1$}}
\put(0,5){\makebox(6,4)[c]{\small \1}}
\put(0,10){\makebox(6,4)[c]{\small \1}}
\put(0,15){\makebox(6,4)[c]{\small \1}}
\put(0,20){\makebox(5.5,5)[c]{\tiny $0|1$}}
\put(5,0){\makebox(6,4)[c]{\small \1}}
\put(5,5){\makebox(6,4)[c]{\small \1}}
\put(5,10){\makebox(6,4)[c]{\small \1}}
\put(5,15){\makebox(5.5,5)[c]{\tiny $0|1$}}
\put(10,0){\makebox(6,4)[c]{\small \1}}
\put(10,5){\makebox(6,4)[c]{\small \1}}
\put(10,10){\makebox(6,4)[c]{\tiny $0|1$}}
\put(15,0){\makebox(6,4)[c]{\small \1}}
\put(15,5){\makebox(6,4)[c]{\tiny $0|1$}}
\put(20,0){\makebox(6,4)[c]{\tiny $0|1$}}
%%%%%%%%%%%%%%%%%%%%%%%%%%%%%%%%%%%%%%%%ETIQUETAGE DES COLONNES ET LIGNES
\put(0,0){\makebox(5,-6)[c]{\tiny $1$}}
\put(5,0){\makebox(5,-6)[c]{\tiny $2$}}
\put(10,0){\makebox(5,-6)[c]{\tiny $3$}}
\put(15,0){\makebox(5,-6)[c]{\tiny $4$}}
\put(20,0){\makebox(5,-6)[c]{\tiny $5$}}
\put(25,0){\makebox(5,-6)[c]{\tiny $6$}}
%%%%%%%%%%%%
\put(0,0){\makebox(-6,5)[c]{\tiny $6$}}
\put(0,5){\makebox(-6,5)[c]{\tiny $5$}}
\put(0,10){\makebox(-5,5)[c]{\tiny $4$}}
\put(0,15){\makebox(-5,5)[c]{\tiny $3$}}
\put(0,20){\makebox(-5,5)[c]{\tiny $2$}}
\put(0,25){\makebox(-5,5)[c]{\tiny $1$}}
\end{picture}}
\hspace{2cm} {\setlength{\unitlength}{1mm}
\begin{picture}(30,35)(0,-5)
%%%%%%%%%%%%%%%%%%%%%%%%%%%%%%%%%%%%%%%%DIAGRAMME
\put(0,0){\line(1,0){30}}\put(0,5){\line(1,0){25}}\put(0,10){\line(1,0){20}}
\put(0,15){\line(1,0){15}}\put(0,20){\line(1,0){10}}\put(0,25){\line(1,0){5}}
\put(0,0){\line(0,1){30}}\put(5,0){\line(0,1){25}}\put(10,0){\line(0,1){20}}
\put(15,0){\line(0,1){15}}\put(20,0){\line(0,1){10}}\put(25,0){\line(0,1){5}}
%%%%%%%%%%%%%%%%%%%%%%%%%%%%%%%%%%%%%%%%REMPLISSAGE
\put(0,0){\makebox(6,4)[c]{\small \1}}
\put(0,5){\makebox(6,4)[c]{\small \1}}
\put(0,10){\makebox(6,4)[c]{\small \1}}
\put(0,15){\makebox(5.5,5)[c]{\tiny $0|1$}}
\put(0,20){\makebox(5.5,5)[c]{\tiny $0|1$}}
\put(5,0){\makebox(6,4)[c]{\small \1}}
\put(5,5){\makebox(6,4)[c]{\small \1}}
\put(5,10){\makebox(6,4)[c]{\small \1}}
\put(5,15){\makebox(6,4)[c]{\small \1}}
\put(10,0){\makebox(6,4)[c]{\small \1}}
\put(10,5){\makebox(6,4)[c]{\small \1}}
\put(10,10){\makebox(6,4)[c]{\small \1}}
\put(15,5){\makebox(6,4)[c]{\small \1}}
\put(15,0){\makebox(5.5,5)[c]{\tiny $0|1$}}
\put(20,0){\makebox(5.5,5)[c]{\tiny $0|1$}}
%%%%%%%%%%%%%%%%%%%%%%%%%%%%%%%%%%%%%%%%ETIQUETAGE DES COLONNES ET LIGNES
\put(0,0){\makebox(5,-6)[c]{\tiny $1$}}
\put(5,0){\makebox(5,-6)[c]{\tiny $2$}}
\put(10,0){\makebox(5,-6)[c]{\tiny $3$}}
\put(15,0){\makebox(5,-6)[c]{\tiny $4$}}
\put(20,0){\makebox(5,-6)[c]{\tiny $5$}}
\put(25,0){\makebox(5,-6)[c]{\tiny $6$}}
%%%%%%%%%%%%
\put(0,0){\makebox(-6,5)[c]{\tiny $6$}}
\put(0,5){\makebox(-6,5)[c]{\tiny $5$}}
\put(0,10){\makebox(-5,5)[c]{\tiny $4$}}
\put(0,15){\makebox(-5,5)[c]{\tiny $3$}}
\put(0,20){\makebox(-5,5)[c]{\tiny $2$}}
\put(0,25){\makebox(-5,5)[c]{\tiny $1$}}
\end{picture}}
\end{center}

\qed

One can also ask if Theorem~\ref{thm:sym-filling}, or more generally
Corollary~\ref{thm:refinement main1} or
Corollary~\ref{thm:refinement main2}, can be extended to arbitrary
larger classes of polyominoes. We note that the condition of
intersection free is necessary. Indeed, the polyomino $T$
represented below is convex but not intersection free,
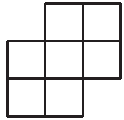
\begin{center}
{\setlength{\unitlength}{1mm}
\begin{picture}(15,20)(0,0)
%%%%%%%%%%%%%%%%%%%%%%%%%%%%%%%%%%%%%%%%DIAGRAMME
\put(0,0){\line(1,0){10}}\put(0,5){\line(1,0){15}}\put(0,10){\line(1,0){15}}
\put(5,15){\line(1,0){10}}
\put(0,0){\line(0,1){10}}\put(5,0){\line(0,1){15}}\put(10,0){\line(0,1){15}}
\put(15,5){\line(0,1){10}}
\end{picture}}
\end{center}
and
$$\sum_{F\in\N^c(T,(1,1,1))}p^{{\nec}_2(F)}q^{{\sec}_2(F)}
=\sum_{F\in\N(T;3)}p^{{\nec}_2(F)}q^{{\sec}_2(F)}=p^2+2q$$ is not
symmetric. One can also check that $\sum_{F\in\N^c(T)}
p^{{\nec}_2(F)}q^{{\sec}_2(F)}$ and $\sum_{F\in\N(T)}
p^{{\nec}_2(F)}q^{{\sec}_2(F)}$ are even not
symmetric.\\

Let $T$ be a moon polyomino and $F$ be a $01$-filling a $T$. Recall
that ${\nec}(F)$ (resp., ${\sec}(F)$) is the largest $k$ for which
$F$ has a NE-chain (resp., SE-chain) of length $k$ such that the
smallest rectangle containing the chain is contained in $F$.
Rubey~\cite{Ru}, answering a conjecture of Jonsson~\cite{Jo}, have
proved that for any positive integers $j$ and $k$, we have
\begin{align}\label{eq:rubey}
&|\{F\in\N^{01}(T) :|F|=j\,,\,{\nec}(F)=k\}|=|\{F\in\N^{01}(T^*)
:|F|=j\,,\,{\nec}(F)=k\}|
\end{align}
for any moon polyomino $T^*$ obtained from $T$ by permutating the
columns (or equivalently the rows) of $T$.

On the other hand, it is easy to derive from
Theorem~\ref{thm:distribution} the following result.
\begin{prop}
Let $T$ be a moon polyomino. For any moon polyomino $T^*$ obtained
from $T$ by permutating the rows of $T$ and any positive integers
$j$, $k$ and $\ell$, we have
\begin{align*} &|\{F\in\N^c(T)
:|F|=j\,,\,{\nec_2}(F)=k\,,\,{{\sec}_2(F)}=\ell\}|\\
=&|\{F\in\N^c(T^*)
:|F|=j\,,\,{\nec_2}(F)=k\,,\,{{\sec}_2(F)}=\ell\}|
\end{align*}
and
\begin{align*} &|\{F\in\N(T)
:|F|=j\,,\,{\nec_2}(F)=k\,,\,{{\sec}_2(F)}=\ell\}|\\
=&|\{F\in\N(T^*) :|F|=j\,,\,{\nec_2}(F)=k\,,\,{{\sec}_2(F)}=\ell\}|.
\end{align*}
\end{prop}

Clearly the above proposition and Rubey's result~\eqref{eq:rubey}
bring us to the following problem: Is it true that for any moon
polyomino $T$ and positive integers $j$ and $k$ we have that
 $$|\{F\in\N^{01}(T) :|F|=j\,,\,{\nec_2}(F)=k\}|
=|\{F\in\N^{01}(T^*) :|F|=j\,,\,{\nec_2}(F)=k\}|$$ for any moon
polyomino $T^*$ obtained from $T$ by permutating the rows of $T$?
The answer is no. Indeed, if such a result holds, then by reflecting
each moon polyomino in a vertical line and apply the result, we
would obtain that the statistics ${\nec_2}$ and ${\nec_2}$ are
equidistributed over ${\N}^{01}(T)$ for any moon polyomino $T$,
which contradicts Proposition~\ref{prop:nonsym-graphe}.

%%%%%%%%%%%%%%%%%%%%%%%%%%%%%%%%%%%%%%%%%%%%%%%%%%%%%%%%%%%%%%%%%%%%%%%%%%%%%%%%%%%%%%%%%%%%%%%%%%%%%%%%%%%%%%%%%%%%%
%%%%%%%%%%%%%%%%%%%%%%%%%%%%%%%%%%%%%%%%%%%%%%%%%%%%%%%%%%%%%%%%%%%%%%%%%%%%%%%%%%%%%%%%%%%%%%%%%%%%%%%%%%%%%%%%%%%%%

\renewcommand{\baselinestretch}{1}


\begin{thebibliography}{99}
\small \setlength{\itemsep}{-.8mm}
%%%%%%%%%%%%%%%%%%%%%%%%%%%%%%%%%%%%%%%%%%%%%%%%%%%%%%%%%%%%%%%%%%%
\bibitem{Bou} M. Bousquet-M\'{e}lou and G. Xin, On partitions avoiding
3-crossings,  S\'{e}m. Lothar. Combin.  54 (2005/07), Art. B54e.
%%%%%%%%%%%%%%%%%%%%%%%%%%%%%%%%%%%%%%%%%%%%%%%%%%%%%%%%%%%%%%%%%%%
\bibitem{Chen1} W. Y. C. Chen, S. Y. J. Wu and C. H. Yan,
linked partitions and linked cycles, European Journal of
combinatorics 29 (2008), Issue 6, 1377-1520.
%%%%%%%%%%%%%%%%%%%%%%%%%%%%%%%%%%%%%%%%%%%%%%%%%%%%%%%%%%%%%%%%%%%
\bibitem{Chen2} W. Y. C. Chen, E. Y. P. Deng, R. R. X. Du, R. P. Stanley and C.
H. Yan, Crossings and nestings of matchings and partitions, Trans.
Amer. Math. Soc. 359 (2007), 1555--1575.
%%%%%%%%%%%%%%%%%%%%%%%%%%%%%%%%%%%%%%%%%%%%%%%%%%%%%%%%%%%%%%%%%%%
\bibitem{Co} S. Corteel,  Crossings and alignments of
permutations, Adv. in Appl. Math.  38  (2007),  no.~2, 149--163.
%%%%%%%%%%%%%%%%%%%%%%%%%%%%%%%%%%%%%%%%%%%%%%%%%%%%%%%%%%%%%%%%%%%
\bibitem{Mi1} A. De Mier, $k$-noncrossing and $k$-nonnesting graphs and fillings
of Ferrers diagrams,  Combinatorica  27  (2007),  no. 6, 699--720.
%%%%%%%%%%%%%%%%%%%%%%%%%%%%%%%%%%%%%%%%%%%%%%%%%%%%%%%%%%%%%%%%%%%
\bibitem{Mi2} A. De Mier, On the symmetry of the distribution of crossings
and nestings in graphs, Electron. J. Combin. 13 (2006), Note \#21.
%%%%%%%%%%%%%%%%%%%%%%%%%%%%%%%%%%%%%%%%%%%%%%%%%%%%%%%%%%%%%%%%%%%
\bibitem{Sai} M. De Sainte-Catherine, Couplage et Pfaffiens en combinatoire,
physique et informatique,
 Th\`ese du 3me cycle, Universit\'e de Bordeaux I, 1983.
%%%%%%%%%%%%%%%%%%%%%%%%%%%%%%%%%%%%%%%%%%%%%%%%%%%%%%%%%%%%%%%%%%%
\bibitem{Jo} J. Jonsson, Generalized triangulations and diagonal-free subsets
of stack polyominoes, Journal of Combinatorial Theory, Series A, 112
(2005), no. 1, 117--142.
%%%%%%%%%%%%%%%%%%%%%%%%%%%%%%%%%%%%%%%%%%%%%%%%%%%%%%%%%%%%%%%%%%%
\bibitem{KaZe} A. Kasraoui and J. Zeng, Distribution of crossings, nestings
and alignments of two edges in matchings and partitions, Electron.
J. Combin. 13 (2006), no. 1, Research Paper 33.
%%%%%%%%%%%%%%%%%%%%%%%%%%%%%%%%%%%%%%%%%%%%%%%%%%%%%%%%%%%%%%%%%%%
\bibitem{KlNo} M. Klazar,  On identities concerning the numbers
of crossings and nestings of two edges in matchings,
  SIAM J. Discrete Math. 20  (2006),  no. 4, 960--976.
%%%%%%%%%%%%%%%%%%%%%%%%%%%%%%%%%%%%%%%%%%%%%%%%%%%%%%%%%%%%%%%%%%%
\bibitem{Kr} C. Krattenthaler, Growth diagrams, and increasing and decreasing
chains in fillings of Ferrers shapes,
  Adv. in Appl. Math. 37  (2006),  no. 3, 404--431.
%%%%%%%%%%%%%%%%%%%%%%%%%%%%%%%%%%%%%%%%%%%%%%%%%%%%%%%%%%%%%%%%%%%
\bibitem{Poz} S. Poznanovik and C. Yan, Crossings and Nestings of Two Edges in Set
Partitions, preprint, avilable at arXiv:0710.1816.
%%%%%%%%%%%%%%%%%%%%%%%%%%%%%%%%%%%%%%%%%%%%%%%%%%%%%%%%%%%%%%%%%%%
\bibitem{Ru} M. Rubey, Increasing and Decreasing Sequences in Fillings of Moon Polyominoes,
preprint, available at arXiv:math.CO/0604140
%%%%%%%%%%%%%%%%%%%%%%%%%%%%%%%%%%%%%%%%%%%%%%%%%%%%%%%%%%%%%%%%%%%

\end{thebibliography}
\end{document}